\documentstyle[11pt]{article}

\newcommand{\ddate}{10 avril 1999}

\newtheorem{th}{Th\'eor\`eme}[section]

\newtheorem{Theorem}[th]{Th\'eor\`eme}

\newtheorem{Lemme}[th]{Lemme}
\newtheorem{Proposition}[th]{Proposition}

\newtheorem{Corollaire}[th]{Corollaire}

\newtheorem{ccote}[th]{}

\newcommand{\preu}{\noindent {\sc Preuve: \ }}
\newcommand{\fl}[1]{\buildrel{#1}\over{\longrightarrow}}

\newcommand{\cqfd}{\unskip\kern 6pt\penalty 500
\raise -2pt\hbox{\vrule\vbox to10pt{\hrule width
 4pt\vfill\hrule}\vrule}\smallskip}

\newcommand{\bbr}{{\bf R}}

\newcommand{\bbz}{{\bf Z}}

\newcommand{\bbn}{{\bf N}}

\newcommand{\cale}{{\cal E}}

\newcommand{\calg}{{\cal G}}
\newcommand{\calh}{{\cal H}}

\newcommand{\calk}{{\cal K}}

\newcommand{\calo}{{\cal O}}

\newcommand{\calr}{{\cal R}}
\newcommand{\cals}{{\cal S}}
\newcommand{\calt}{{\cal T}}

\newcommand{\calv}{{\cal V}}

\newcommand{\bft}{{\bf T}}

\newcommand{\bfw}{{\bf W}}

\newcommand{\pcirc}{\kern .7pt {\scriptstyle \circ} \kern 1pt}
\newcommand{\iso}{\cong}

\newcommand{\hfl}[1]{\buildrel{#1}\over{\longrightarrow}}
\newcommand{\eqref}[1]{(\ref{#1})}
\newcounter{exo}

\newcommand{\gro}{groupo\"\i de }
\newcommand{\gros}{groupo\"\i des }
\newcommand{\cb}{$\clubsuit$}

\newcommand{\SIN}{\hbox{\sl SIN}}
\newcommand{\res}{\hbox{\rm res\,}}
\newcommand{\calrr}{\calr(\Omega_C,G)_\xi}
\newcommand{\calar}{\calr(C,\xi)}

\newcommand{\mapr}{{\rm Map\,}^{\scriptscriptstyle\bullet}(X,BG)_\xi}
\newcommand{\maper}{{\rm Map\,}_G^{\scriptscriptstyle\bullet}(E,EG)}
\newcommand{\maprs}{{\rm Map\,}^{\scriptscriptstyle\bullet}(X,B)_\xi}
\newcommand{\mapers}{{\rm Map\,}_G^{\scriptscriptstyle\bullet}(E,A)}
\newcommand{\llangle}[2]{\langle #1 ,#2 \rangle}

\title{Th\'eorie de jauge et groupo\"\i des}
\author{Jean-Claude HAUSMANN}
\date{\ddate}

\begin{document}    \maketitle
\begin{abstract}
Les r\'esultats de cet article concernent le probl\`eme de l'existence 
de repr\'esentations d'un \gro\ topologique sur un fibr\'e principal et leur classification
\`a transformation de jauge pr\`es. De telles repr\'esentations interviennent naturellement dans 
divers contextes (th\'eo\-ries de jauge 
classiques ou sur graphe, fibr\'es \'equivariants, etc).
\end{abstract}

\tableofcontents

\section{Pr\'esentation des r\'esultats} \label{intro}

Soit $X$ un espace topologique. 
Un $X$-{\it \gro\ } est un espace topologique
$C$ muni de :

a) deux applications continues $\alpha,\beta : C\to X$ ({\it source} et {\it but}).
   Pour $x,y\in X$, on note $C_x:=\alpha^{-1}(\{x\})$, $C^y:=\beta^{-1}(\{y\})$
   et $C_x^y:= C_x\cap C^y$. Le \gro\ $C$ est dit {\it transitif} si 
   $C_x^y\not=\emptyset$ pour tout $x,y\in X$.

b) une {\it composition} partiellement d\'efinie
$$
C\times_X C :=\{(c_1,c_2)\in C\mid \alpha(c_1)=\beta (c_2)\}  \to  C
\quad ,\quad (c_1,c_2) \mapsto  c_1c_2 \kern 3 truecm 
$$
qui est continue, $C\times_X C$ \'etant un sous-espace de $C\times C$ avec la topologie
produit. On demande que $\alpha(c_1c_2)=\alpha(c_2)$, $\beta(c_1c_2)=\beta(c_1)$ et que la
composition soit associative. 

c) une application continue $i : X\to C$, associant \`a $x$ l'{\it unit\'e}
$i_x\in C^x_x$ qui est \'el\'ement neutre \`a gauche et \`a droite pour la composition.

d) une anti-involution continue $c\mapsto c^{-1}$ de $C$, 
envoyant $C^y_x$ sur $C^x_y$, telle que $cc^{-1}=i_{\alpha(c)}$ et 
$c^{-1}c=i_{\beta(c)}$. 

On suppose que $X$ est muni d'un point base $*\in X$. En vertu de d), l'espace
$C^*_*$ est un groupe topologique que l'on note $\Omega_C$.

\vskip .2 truecm

On peut voir $C$ comme l'espace des morphismes (tous inversibles)
d'une petite cat\'egorie topologique dont l'espace des objets est $X$. 
C'est le point de vue de \cite{Ma}. 
Si l'on pr\'ef\`ere les \gros\ "sans objets" \cite[II.5]{Co},
on pr\'esentera $C$ comme un 
\gro\ topologique dont l'espace des unit\'es est identifi\'e \`a $X$.

Soit $G$ un groupe topologique et
soit $\xi:= (E\fl{p}{} X)$ un {\it $G$-espace principal au dessus de $X$}.
On entend par l\`a que $p$ est continue et que l'on s'est donn\'e une action continue libre 
$E\times G\to E$ telle que $p$ induise une bijection continue 
du quotient $E/G$ sur $X$. 
Par exemple, un $G$-fibr\'e principal sur $X$ 
est un $G$-espace au dessus de $X$ qui est localement trivial.
On suppose que l'espace $E$ est point\'e par $\tilde * \in p^{-1}(*)$. 

Une {\it repr\'esentation} d'un $X$-\gro\ $C$ 
 sur un $G$-espace $\xi$ est 
une application continue $w: C\times_X E\to E$ 
(o\`u $C$ est vu au dessus de $X$ via $\alpha$) telle que, pour tout 
$c,d\in C$, $z\in E$ et $g\in G$, on ait
\begin{enumerate}
\item $p(w(c,z)) = \beta(c)$.
\item $w(c\, d,z)=w(c,w(d,z))$.
\item $w(c^{-1},(w(c,z))=z$.
\item $w(c,z\cdot g) = w(c,z)\cdot g$.
\end{enumerate}

Par exemple, si $X=*$, alors $C$ est un groupe topologique et $w$ correspond \`a
une repr\'esentation (= homomorphisme continu) de $C$ dans $G$. Plus g\'en\'eralement, 
si $w$ est une repr\'esentation de $C$ sur $\xi$,
l'\'equation $w(c,\tilde *) = \tilde * \cdot h(c)$ d\'efinit un
homomorphisme continu $h^w:\Omega_C\to G$ appel\'e l'{\it holonomie} de $w$.

Les r\'esultats de cet article concernent 
l'existence et la classification de repr\'esentations d'un $X$-\gro\ $C$ 
 sur un $G$-fibr\'e principal
$\xi$ donn\'es. Nos hypoth\`eses principales seront que $C$ est un $X$-\gro\ 
{\it localement trivial} et $G$ un groupe \SIN\ (d\'efinitions ci-dessous).

Soit $C$ un $X$-groupo\"\i de. Une {\it $C$-contraction}
 est une application continue $\rho : U_\rho\to C^*$, o\`u  
$U_\rho$ est un ouvert de $X$, telle que $\alpha(\rho(x))=x$. Soit
${\rm Cont\,}_C(X)$ l'ensemble des $C$-contractions. 
Le $X$-\gro $C$ est dit {\it localement trivial} si
$\{U_\rho\mid \rho\in{\rm Cont\,}_C(X)\}$ est un recouvrement
ouvert de $X$ (\cite[p. 32]{Ma}; cette d\'efinition co\"\i ncide avec la terminologie de 
\cite{Eh} si $C$ est transitif).
Si le recouvrement $\{U_\rho\}$ est de plus 
num\'erisable\footnote{Nous utilisons ``num\'erisable" pour \'equivalent fran\c cais 
du n\'eologisme anglais ``numerable" introduit par Dold [Do].
De m\^eme, un fibr\'e sera dit {\it num\'erisable} s'il admet un recouvrement 
num\'erisable d'ouverts trivialisants.}, 
c'est-\`a-dire s'il existe une partition de l'unit\'e 
$\{\mu_\rho : X\to [0,1]\}$ qui lui est subordonn\'ee, on dira 
que $C$ est un {\it $X$-\gro\ localement trivial num\'erisable}.
Plusieurs exemples naturels de  tels \gros sont pr\'esent\'es
au paragraphe \ref{P:exp}.

Reprenant les id\'ees de Ch. Ehresmann \cite{Eh}, 
nous d\'emontrerons au \S\ \ref{P:pthma}
l'existence d'une {\it repr\'esentation universelle de $C$}
et d'un principe de reconstruction~:

\begin{Theorem}[Th\'eor\`eme A] \label{couniv}
Soit $C$ un $X$-\gro localement trivial num\'erisable. Alors~:

a) l'application $\beta: C_*\to X$ est un $\Omega_C$-fibr\'e principal
num\'erisable (que nous appellerons $\xi_C$). 
La formule $\tilde w(c,u):=c\, u$ d\'efinit une repr\'esentation de $C$ sur $\xi_C$.

b) soit $\xi:= (E\fl{p}{} X)$ un $G$-espace principal muni d'une
repr\'esentation  $w$ de $C$. Alors $\xi$ est le fibr\'e associ\'e \`a $\xi_C$ 
$$E= C_* \times _{\Omega_C} G = C_* \times _{h^w} G \, ,$$
o\`u  $\Omega_C$ agit \`a gauche sur $G$ 
via l'holonomie $h^w$ de $w$.
En particulier, $\xi$ est un $G$-fibr\'e principal num\'erisable.
De plus, la repr\'esentation $w$ est obtenue de $\tilde w$ par la formule 
$w(c,[u,g])=[\tilde w(c,u),g]=[cu,g]$.
\end{Theorem}

Rappelons qu'un $G$-fibr\'e principal num\'erisable $\xi$ sur $X$
est induit du fibr\'e universel $EG\to BG$ sur le classifiant de Milnor
$BG$ par une application continue (unique \`a homotopie pr\`es) $\nu_\xi :X\to BG$.
En particulier, d'apr\`es le th\'eor\`eme A, un $X$-\gro localement trivial num\'erisable $C$
donnera une application $\nu_C:=\nu_{\xi_C}:X\to B\Omega_C$ 
induisant $\xi_C$. Le th\'eor\`eme d'existence d'une repr\'esentation de $C$ sur $\xi$,
qui d\'ecoule facilement du th\'eor\`eme A, est le suivant~:

\begin{Theorem}[Th\'eor\`eme d'existence] \label{counivex}
Soit $C$ un $X$-\gro localement trivial num\'erisable
et $\xi$ un $G$-fibr\'e principal.
Alors, $\xi$ admet une repr\'esentation de $C$ si et seulement si
il est num\'erisable et s'il existe un homomorphisme continu
$\phi : \Omega_C\to G$ tel que $B\phi\pcirc \nu_C$ soit homotope \`a $\nu_\xi$.
\end{Theorem}

Ayant r\'esolu la question de l'existence, int\'eressons-nous \`a l'ensemble
$\calar$ des repr\'esentations de $C$ sur $\xi$. Il est muni d'une action
du {\it groupe de jauge} $\calg$ de $\xi$. Les \'el\'ements de $\calg$,
les {\it transformations de jauge}, sont les automorphismes 
$G$-\'equi\-va\-riants de $E$ au dessus de ${\rm id}_X$.
Pour $w\in\calar$ et $\chi\in\calg$, on d\'efinit $w^\chi$ par
$$w^\chi(c,z):=\chi^{-1}(w(c,\chi(z)).$$
ce qui donne une action \`a droite de $\calg$ sur $\calar$.
Le sous-groupe invariant $\calg_1$ de $\calg$ 
form\'e des transformations de jauge $\tilde *$ 
joue un r\^ole important car il 
agit librement sur $\calar$ (voir \ref{L:g1lib}). 

Les ensembles $\calar$ et $\calg$ sont munis de la topologie
compact-ouvert (CO-topologie).
Nous d\'emontrerons, au \S\ \ref{unif} que $\calg$ est un groupe
topologique et que l'action $\calar\times\calg\to\calar$ est continue,
ceci sous l'hypoth\`ese que $G$ est \SIN . 
Rappelons qu'un groupe topologique est 
\SIN , si son \'el\'ement neutre admet un syst\`eme fondamental de voisinages
qui sont $G$-invariants (par conjugaison; \SIN\ = small invariant neighbourhood).
Par exemple, les groupes compacts, ab\'elien ou discrets sont \SIN\
(voir \cite{Pa}, pour la litt\'erature classique sur ces groupes).
Nous montrerons que si $G$ est \SIN, la CO-topologie sur $\calar$ et $\calg$
provient de structures uniformes et que toutes les op\'erations
sont uniform\'ement continue (voir \S\ \ref{unif}).

Pour \'etudier les quotients $\calar/\calg_1$ et $\calar/\calg$, 
introduisons l'espace $\calr (\Omega_C,G)$ des
repr\'esentations (i.e. homomorphismes continus) de $\Omega_C$ dans $G$, 
muni de la CO-topologie.
Soit $\calr (\Omega_C,G)_\xi$ le sous-espace des $\phi\in\calr (\Omega_C,G)$
tels que l'application compos\'ee $n(\phi):=B\phi\pcirc \nu_C$ soit homotope \`a $\nu_\xi$
($\calr (\Omega_C,G)_\xi$ est une union de componates connexes 
par arc de $\calr (\Omega_C,G)$). 

Au niveau des espaces classifiants, on d\'esigne
 par $\mapr$ l'espace des applications
continues point\'ees de $f:X\to BG$ 
qui induisent $\xi$, muni de la CO-topologie.

Un espace topologique $Z$ est dit {\it semi-localement contractile}
si tout $z\in Z$ admet un voisinage $U_z$ dont l'inclusion $U_z\hookrightarrow Z$ est
homotope \`a une application constante. Cette condition ne d\'epend que du type 
d'homotopie de $Z$.

\begin{Theorem}[Th\'eor\`eme B] \label{gau}
Soit $C$ un $X$-\gro localement trivial num\'erisable et s\'epar\'e. Soit $\xi$ 
un $G$-fibr\'e principal num\'erisable sur $X$.
On suppose que $G$ est \SIN , que $X$ est
localement compact et que 
$\mapr$ est semi-localement contractile.
Soit  $h:\calar\to\calr (\Omega_C,G)_\xi$ 
l'application qui \`a une repr\'esentation de $C$ 
$w$ associe son holonomie $h^w$. Alors $h$ est
est un $\calg_1$-fibr\'e principal.
\end{Theorem}

Nous ne savons pas si l'hypoth\`ese 
``$\mapr$ semi-localement contractile" est toujours v\'erifi\'ee. 
Observons que cette condition 
ne d\'epend que des types d'homotopie de 
$X$ et de $G$. Elle est vraie si, par exemple, $X$ est compact et 
$G$ est un groupe de Lie compact. En effet, $BG$ a alors le type d'homotopie 
d'un CW-complexe d\'enombrable, limite inductive quotients 
de vari\'et\'e de Stiefel \cite[\S\ 19.6]{St}.
L'espace $\mapr$ est ainsi semi-localement contractile par 
\cite[lemmes 2 et 3 p. 277]{Mi1}. D'autre part, 
tout recouvrement d'un espace compact est num\'erisable, ce qui prouve le~:

\begin{Corollaire}\label{gaucor}
Soit $C$ un $X$-\gro localement trivial et s\'epar\'e, avec $X$ compact. Soit $\xi$ 
un $G$-fibr\'e principal sur $X$ avec $G$ un groupe de Lie compact.
Alors  $h:\calar\to\calr (\Omega_C,G)_\xi$ est
est un $\calg_1$-fibr\'e principal.
\end{Corollaire}

Le th\'eor\`eme B montre que le quotient $\calar/\calg_1$ est hom\'eomorphe \`a
$\calr (\Omega_C,G)_\xi$. Pour d\'ecrire  
$\calar/\calg$, on consid\`ere l'action \`a droite de $G$ sur 
$\calr (\Omega_C,G)$ par conjugaison. 
Si $G$ est connexe par arc, cette action 
pr\'eserve $\calr (\Omega_C,G)_\xi$.

\begin{Theorem}[Th\'eor\`eme C] \label{gauser}
Avec les hypoth\`eses du th\'eor\`eme B, l'appli\-ca\-tion
$h$ induit un hom\'eomorphisme
$\calar /\calg_1\hfl{\iso}{} \calr (\Omega_C,G)_\xi$. 
Si, de plus, $G$ est connexe par arc, elle induit un hom\'eomorphisme
$\calar /\calg\hfl{\iso}{} \calr (\Omega_C,G)_\xi/G.$ 
\end{Theorem}

Les th\'eor\`emes B et C pr\'esentent une analogie avec des r\'esultats 
de la th\'eorie de jauge qu'il serait int\'eressant d'\'etudier plus \`a fond.
On sait que, si $\xi$ est un 
$G$-fibr\'e diff\'erentiable sur une vari\'et\'e compacte $X$ ($G$ groupe de Lie
compact connexe), alors l'espace des connexions sur $\xi$ modulo 
$\calg_1$ a le type d'homotopie faible de $\mapr$ \cite[Prop. 5.1.4]{DK}. 
On verra au \S\ \ref{clipamconn} qu'une connexion donne 
une repr\'esentation d'un \gro\ ${\bf D}(X)$ 
construit \`a l'aide des chemins  dans $X$ lisses par morceau. 
Pour l'instant, observons que les espaces $\calrr$ et $\mapr$ sont reli\'es par une application continue
$n:\calrr\to\mapr$ (d\'etails en \ref{applicn}) dont on peut d\'ecrire 
la fibre homotopique lorsque $X$ est ``bien point\'e"~:

\begin{Theorem}[Th\'eor\`eme D] \label{gauserthd}
Supposons que l'on ait les hypoth\`eses du th\'eor\`eme B
et que l'inclusion $\{*\}\subset X$ soit une cofibration. Alors,
la fibre homotopique de $n : \calr (\Omega_C,G)_\xi \to \mapr$ 
a le type d'homotopie faible de $\calar$.
\end{Theorem}

Le \S\ \ref{P:pthma} est consacr\'e \`a la preuve du th\'eor\`eme A et du th\'eor\`eme d'existence.
Les \S\ \ref{unif} et \ref{P:gfibuniv} pr\'eparent aux preuves des th\'eor\`emes B,C et D 
qui sont donn\'ees dans le \S\ \ref{P:pthmbc}. Enfin, le \S\ \ref{P:exp} pr\'esente quelques
exemples et applications.

\paragraph{Remerciements : } Ce travail a b\'en\'efici\'e du support du Fonds National 
Suisse de la Recherche Scientifique.
L'auteur tient \'egalement \`a remercier E. Dror-Farjoun, P. de
la Harpe et R. Vogt pour d'utiles discussions.

\section{Preuve des th\'eor\`emes A et d'existence} \label{P:pthma}

\begin{Lemme} \label{rectrican}
Soit $C$ un $X$-\gro localement trivial num\'erisable.
Soit $\xi:(E\hfl{p}{}X)$ un $G$-espace principal 
admettant une repr\'esentation $w$ de $C$. Alors, 
$\xi$ est un $G$-fibr\'e principal avec recouvrement trivialisant
(num\'erisable) $\{U_\rho\mid\rho\in
{\rm Cont\,}_C(X)\}$.  Plus pr\'ecis\'ement, 
$\xi$ restreint \`a $U_\rho$ admet la trivialisation 
$\psi^w_\rho: p^{-1}(U_\rho)\to U_\rho\times G$ donn\'ee par
\begin{equation}\label{deftriv1}
\psi^w_\rho(z) := (p(z), \lambda_w(z)) \quad \hbox{ o\`u } \quad
w(\rho(p(z)),z) = \tilde * \lambda_w(z).\end{equation}
\end{Lemme}

\preu Il est banal que $\psi^w_\rho$ est $G$-equivariante et on
 v\'erifie que inverse de $\psi^w_\rho$ est 
$(\psi^w_\rho)^{-1}(x,g):=w(\rho(x)^{-1}, \tilde * )\cdot g$  \cqfd.

\paragraph{Preuve du th\'eor\`eme A}
Observons que l'application $\beta:C_*\to X$ est surjective si 
$C$ est localement trivial. En effet, si $x\in X$, il existe $\rho\in  {\rm
Cont\,}_C(X)$ telle que $x\in U_\rho$ et alors $x=\beta(\rho(x)^{-1})$.
L'action de $\Omega_C$ sur $C_*$ est d\'efinie par
$u\cdot c := u\, c$. Si $\beta(u)=\beta(\bar u)$, alors 
$\bar u = \bar u\cdot c$ avec $c:= (\bar u^{-1}\, u)$. L'application
$\beta$ induit donc une bijection continue  
$\bar \beta :C_*/\Omega_C\to X$. Les trivialistions construites 
ci-dessous font que $\beta$ est ouverte 
et donc $\bar\beta$ est un hom\'eomorphisme.

Pour  $\rho\in {\rm Cont\,}_C(X)$, on d\'efinit la trivialisation
$\hat \psi_\rho : \beta^{-1}(U_\rho)\to U_\rho\times \Omega_C$
de la mani\`ere suivante:
\begin{equation} \label{deftriv2}
\hat\psi_\rho(u) = (\beta(u), \rho(\beta (u))\, u).
\end{equation}
Elle admet pour inverse 
\begin{equation}
\hat\psi^{-1}_\rho(x,c) = \rho(x)^{-1}c.
\end{equation}
La $\Omega_C$-equivariance de l'hom\'eomorphisme $\hat\psi_\rho$ est banale.
Ceci montre que $\beta : C_*\to X$ est un fibr\'e $\Omega_C$-principal avec
recouvrement trivialisant num\'erisable $\{U_\rho\}$. 
Le fait que la formule $\tilde w (c,u)= c\, u$ d\'efinisse une 
repr\'esentation de $C$ sur ce fibr\'e provient des axi\^omes de \gro .

Passons au point b) du th\'eor\`eme A. 
Soit 
$w$ une repr\'esentation de $C$ sur un $G$-espace principal 
$\xi:= (E\fl{p}{} X)$.
Le Lemme \ref{rectrican} montre que
$\xi$ est un fibr\'e principal num\'erisable.
L'espace $E$ \'etant point\'e par $\tilde * \in p^{-1}(*)$,
on d\'efinit $\Phi : C_* \times _{\Omega_C} G \to E$ par
\begin{equation}
\Phi(u,g) = w(u,\tilde * ) \cdot g
\end{equation}

a) {\it $\Phi$ est bien d\'efinie : } Soit $c\in \Omega_C$. On a
$$\begin{array}{lcl}
\Phi(u\, c,g) & = &  w(u\, c,\tilde * ) \cdot g =
w(u,w(c,\tilde * )) \cdot g = w(u, \tilde *\cdot h(c))\cdot g
= \\[2pt] & = & w(u, \tilde *)\cdot (h(c) g) = \Phi(u, h(c)g).\end{array}$$

\vskip .2 truecm
b) {\it $\Phi$ est $G$-equivariante : } \'evident.

\vskip .2 truecm
c) {\it $\Phi$ est surjective : } Soit $z\in E$. 
Soit $\rho\in {\rm Cont\,}_C(X)$ telle que $p(z)\in U_\rho$. Posons
$u:=\rho(p(z))^{-1}\in C_*^{p(z)}$. On a ainsi $p(\Phi(u,1))= \beta(u)=p(z)$. Il existe
donc un unique $g\in G$ tel que $\Phi(u,g)=\Phi(u,1)\cdot g = z$.

\vskip .2 truecm
d) {\it $\Phi$ est injective : } 
Supposons que $\Phi(u,g)=\Phi(\bar u,\bar g)$. On a donc 
$\bar u =  u \, c$ avec $c:=u^{-1}\bar u\in\Omega_C$ et
$$\begin{array}{lcl}
\Phi(\bar u ,\bar g) & = &  w( u\, c,\tilde * ) \cdot \bar g =
w(u,w(c,\tilde * )) \cdot \bar g = 
\\[2pt] & = & w(u, \tilde *)\cdot (h(c) \bar g) 
\end{array}$$
Comme d'autre part
$$\Phi(\bar u ,\bar g)=\Phi(u ,g) = w(u, \tilde *)\cdot g,$$
il s'en suit que $g=h(c) \bar g$. Dans $C_*\times _{\Omega_C} G$, on aura
ainsi les \'egalit\'es
$$(\bar u,\bar g) = (u\, c,\bar g) = (u, h(c)\bar g) = (u,g).$$

\vskip .2 truecm
e) {\it $\Phi$ est un hom\'eomorphisme : }  Soit $\rho\in {\rm Cont\,}_C(X)\}$.
Avec les trivialisations $\psi_\rho$ et $\hat \psi_\rho$
introduites en \eqref{deftriv1} et \eqref{deftriv2}, l'application
$$\Phi\rho:=\psi_\rho \pcirc \Phi \pcirc \hat \psi_\rho^{-1} :
U_\rho\times(\Omega_C\times_{\Omega_C}G)\to U_\rho\times G$$
 s'\'ecrit
$\Phi\rho(x,(c,g)) = (x,h(c)\, g)$. L'application $\Phi_\rho$ est donc
 un hom\'eo\-mor\-phsime pour tout $\rho$, ce qui
implique que $\Phi$ est un hom\'eomorphisme.

\paragraph{Preuve du th\'eor\`eme d'existence}  
Supposons que $\xi$ admette une repr\'e\-sen\-tation $w$ de $C$. 
Par le point b) du th\'eor\`eme A, $\xi$ est obtenu du fibr\'e $\xi_C$ par 
la construction de Borel avec l'holonomie $h^w$. Cela implique
que $\nu_\xi$ est homotope \`a $Bh^w\pcirc \nu_C$.

R\'eciproquement, supposons qu'il existe un homomorphisme continu
$\phi : \Omega_C\to G$ tel que $\nu_\xi$ est homotope \`a
$B\phi\pcirc \nu_C$. L'espace total du fibr\'e induit 
$\nu_xi^*(EG)$ est de la forme 
$C_*\times_{\Omega_C} G$ et admet donc la 
repr\'esentation $w(c,[u,g])=[cu,g]$. Le fibr\'e $\nu_xi^*(EG)$ \'etant
isomorphe \`a $\xi$, cela donne une repr\'esentation de $C$ sur $\xi$. \cqfd

\section{
Structures uniformes sur $\calg$ et $\calar$} \label{unif}

Dans ce paragraphe, $\xi: E\hfl{p}{}X$ d\'esigne un $G$-fibr\'e principal
num\'erisable et $\calg$ son groupe de jauge, muni de la CO-topologie. 
Le fait que $\calg$ est un groupe topologique 
est non-trivial, car $E$ n'est m\^eme pas suppos\'e localement compact.
Pour d\'emontrer ce fait, nous allons, lorsque $G$ est \SIN ,
munir $\calg$ d'une
structure uniforme ${\bf U}_\calg$, dont on montrera, si $X$ est s\'epar\'e,
qu'elle induit
la CO-tolopogie  (Proposition \ref{compou}). La m\^eme strat\'egie sera utilis\'ee
pour l'action de $\calg$ sur $\calar$.

Pour d\'ecrire ${\bf U}_\calg$,
consid\'erons l'application continue $\gamma :E\times_X E\to G$ 
d\'efinie par l'\'equation
\begin{equation}\label{defgamma}
y=z\cdot \gamma(y,z) \quad, \quad (y,z)\in E\times_X E.
\end{equation}
Soit $\calv_G$ l'ensemble des ouverts $V$ 
de $G$ contenant l'\'el\'ement neutre et tels que $V=V^{-1}$.
Pour $V\in\calv_G$ et $K$ un compact de $X$, 
on d\'efinit 
$$\calo^\calg({K,V}) :=\{(\chi,\tilde\chi)\in\calg\times\calg\mid
\gamma(\chi(z),\tilde\chi(z))\in V\hbox{ \small pour tout } 
z\in p^{-1}(K)\}.$$
Comme $\calo^\calg({K_1,V_1})\cap\calo^\calg({K_2,V_2})$ contient
$\calo^\calg({K_1\cup K_2,V_1\cap V_2})$, la famille 
$\calo^\calg(K,V)$ est un syst\`eme fondamental d'entou\-rages
de la diagonale dans $\calg\times\calg$, d\'eterminant, par d\'efinition,
la structure uniforme ${\bf U}_\calg$ sur $\calg$.

\begin{Proposition} \label{unicalg}
Si $G$ est un groupe {\it SIN}, le groupe de jauge $\calg$, 
muni de la structure uniforme ${\bf U}_\calg$, 
est un groupe topologique. De plus, $\calg$ est alors {\it SIN}. 
\end{Proposition}

\preu 
Nous allons tout d'abord d\'emontrer que la multiplication
$\calg\times\calg\to\calg$ 
et le passage \`a l'inverse sont uniform\'ement continus.
La topologie induite par ${\bf U}_\calg$ fera donc de $\calg$ un
groupe topologique dont on v\'erifiera directement qu'il est \SIN .

Soient $\chi_1$, $\chi_2$, $\tilde \chi_1$, $\tilde \chi_2$ des \'el\'ements de
$\calg$.  
Par d\'efinition de l'application $\gamma :E\times_X E\to G$, on a,
pour $z\in E$~:
$$\tilde\chi_1\pcirc\tilde\chi_2(z) =
\chi_1\pcirc\chi_2(z) \cdot
\gamma(\tilde\chi_1\pcirc\tilde\chi_2(z),\chi_1\pcirc\chi_2(z))
.$$
Par $G$-equivariance des transformations de jauge, on a~:
$$\begin{array}{lcl}
\tilde\chi_1\pcirc\tilde\chi_2(z) &=& 
\tilde\chi_1\big(\chi_2(z)\cdot\gamma(\tilde\chi_2(z),\chi_2(z))\big)=
\tilde\chi_1(\chi_2(z))\cdot\gamma(\tilde\chi_2(z),\chi_2(z)) =
\\[2pt] &=&
\chi_1(\chi_2(z))\cdot
\gamma(\tilde\chi_1(\chi_2(z)),\chi_1(\chi_2(z)))
\cdot\gamma(\tilde\chi_2(z),\chi_2(z))
\end{array}.$$
On en d\'eduit que
$$\gamma(\tilde\chi_1\pcirc\tilde\chi_2(z),\chi_1\pcirc\chi_2(z))=
\gamma(\tilde\chi_1\pcirc\chi_2(z),\chi_1\pcirc\chi_2(z))
\,\, \gamma(\tilde\chi_2(z),\chi_2(z)).$$
Soit $V\in\calv_G$. Comme $G$ est un groupe topologique, il existe
$W\in\calv_G$ tel que $W\cdot W \subset V$. 
La condition 
$(\tilde\chi_1\pcirc\tilde\chi_2,\chi_1\pcirc\chi_2)\in\calo^\calg(K,V)$
sera vraie si
$(\tilde \chi_1,\chi_1)$ et $(\tilde \chi_2,\chi_2)$ sont dans $\calo^\calg(K,W)$. 
Ceci d\'emontre la continuit\'e uniforme de la composition 
$\calg\times\calg\to\calg$.

Pour le passage \`a l'inverse, soient $\chi,\tilde\chi\in\calg$.
Les \'equations
$$\chi(z)=z\cdot \gamma(\chi(z),z) \quad , \quad 
\tilde\chi(z)=z\cdot \gamma(\tilde\chi(z),z)$$
donnent
$$\tilde\chi(z)=\chi(z)\cdot \gamma(\chi(z),z)^{-1}
\cdot \gamma(\tilde\chi(z),z)$$
d'o\`u
\begin{equation}\label{E:inv1}
\gamma(\tilde\chi(z),\chi(z))=
\gamma(\chi(z),z)^{-1}\,  \gamma(\tilde\chi(z),z).
\end{equation}
Observons que
$$\chi_1(\chi_2(z))=\chi_1(z)\cdot\gamma(\chi_2(z),z)=
z\cdot\gamma(\chi_1(z),z)\cdot\gamma(\chi_2(z),z)$$
et donc
\begin{equation}\label{E:inv2}
\gamma(\chi_1\pcirc\chi_2(z),z)=
\gamma(\chi_1(z),z)\,\gamma(\chi_2(z),z).
\end{equation}
On en d\'eduit que $\gamma(\chi^{-1}(z),z)=\gamma(\chi(z),z)^{-1}$.
En changeant $\chi,\tilde\chi$ en 
$\chi^{-1},\tilde\chi^{-1}$ dans \eqref{E:inv1}, on obtient ainsi
\begin{equation}\label{E:inv3}
\begin{array}{lcl}
\gamma(\tilde\chi^{-1}(z),\chi^{-1}(z)) &=&
\gamma(\chi(z),z)\,  \gamma(\tilde\chi(z),z)^{-1}=\\[2pt] &=&
\gamma(\chi(z),z)\,\gamma(\tilde\chi(z),\chi(z))^{-1}\,\gamma(\chi(z),z)^{-1}.
\end{array}\end{equation}
Soit $V\in\calv_G$ et $K$ un compact de $X$. Comme $G$ est \SIN , il existe un voisinage invariant
$W$ de l'\'el\'ement neutre contenu dans $V$. En rempla\c cant au besoin 
$W$ par $W\cup W^{-1}$, on peut supposer que $W=W^{-1}$. 
Gr\^ace \`a l'\'equation
\eqref{E:inv3}, si $(\tilde\chi,\chi)\in\calo^\calg(K,W)$, alors
$(\tilde\chi^{-1},\chi^{-1})\in\calo^\calg(K,V)$, 
ce qui prouve la continuit\'e uniforme de $\chi\mapsto\chi^{-1}$.

Pour voir que $\calg$ est \SIN , on utilise que tout voisinage de 
${\rm id}_E$ contient un voisinage du type 
$\calh(K,W):=\{\chi\mid (\chi,{\rm id}_E)\in\calo^\calg(K,W)\}$
o\`u $K$ est un compact de $X$ et $W$ un voisinage 
invariant de l'\'el\'ement neutre dans $G$. Par la formule \eqref{E:inv2},
on a, pour tout $\tilde\chi\in\calg$ et tout $z\in E$
$$\gamma(\tilde\chi^{-1}\pcirc\chi\pcirc\tilde\chi(z),z)=
\gamma(\tilde\chi(z),z)^{-1}\,\gamma(\chi(z),z)\,
\gamma(\tilde\chi(z),z).$$
On en d\'eduit imm\'ediatement que $\calh(K,W)$ est $\calg$-invariant. \cqfd

D\'efinissons $\res : \calg\to G$ par la formule
\begin{equation}\label{E:defres}
     \chi(\tilde *) = \tilde *\cdot\res (\chi).
\end{equation}
L'\'equation \ref{E:inv2} implique que $\res$ est un homomorphisme.
Le noyau de $\res$ est \'evidemment $\calg_1$.
Le groupe $G$ est muni de sa structure uniforme naturelle~: une base
d'entourages est donn\'ee par 
$\calo^V:=\{(\tilde g ,g)\mid \tilde g g^{-1}\in V\}$. On obtient
la m\^eme structure uniforme avec la condition $\tilde g^{-1}g\in V$
lorsque $G$ est \SIN .

\begin{Proposition} \label{numer}
Si $G$ est \SIN , l'homomrophisme $\res$ est uniform\'ement continu.
Si $G$ est connexe par arc, $\res$ est surjectif.
\end{Proposition}

\preu Pour tout $V\in\calv_G$, on a $\res(\calo^\calg(\{*\},V)\subset V$,
ce qui prouve la continuit\'e uniforme de $\res$ (ceci ne semble pas
utiliser que $G$ est \SIN\ mais rappelons que cette condition est
n\'ecessaire pour que $G$ soit un groupe topologique par la proposition
\ref{unicalg}).

La preuve de la surjectivit\'e de $\res$ utilise que $\xi$ est num\'erisable.
Soit $\mu :X\to\bbr$ une application continue telle que 
$\mu(*)\not = 0$ et dont le support est contenu dans un ouvert trivialisant
$U$. Soit $\psi : p^{-1}(U)\to U\times G$ une trivialisation. 
Posons $\psi(z)=(p(z),\delta(z))$. En divisant 
$\mu$ par $\mu(*)$, on peut supposer que $\mu(*)=1$. Soit $g\in G$. Comme
$G$ est connexe par arc, il existe un chemin continu $g(t)$ avec
$g(0)=e$ et $g(1)=g$. On d\'efinit alors $\chi\in\calg$ par
$$\chi(z):=\left\{\begin{array}{lllll}
\psi^{-1}(p(z),g(\mu(p(z))\,\delta(z))\ & \hbox{si }  p(z)\in U\\
z & \hbox{sinon. }\end{array}\right.$$
Il est clair que $\res (\chi) = g$.  \cqfd

Nous allons maintenant munir l'espace $\calar$ des repr\'esentations de $C$ sur $\xi$
de la structure uniforme ${\bf U}_\calar$ dont
une base d'entourages est donn\'ee par
$$\calo^\calr(L,V) :=\{(w,\tilde w)\in\calar\times\calar\mid
\gamma(w(c,z),\tilde w(c,z))\in V
\ \forall (c,z)\in L\times_X E\}$$
o\`u $V\in \calv_G$ et $L$ est un un compact de $C$.

\begin{Proposition} \label{unicala}
L'action $\calar\times\calg\to\calar$ est uniform\'ement continue.
\end{Proposition}

\preu Soient $\chi,\tilde\chi\in\calg$ et $w,\tilde w\in\calar$. 
Pour $(c,z)\in C\times_X E$, on d\'emontre, comme dans la preuve de la proposition
\ref{unicala} que
$\gamma(\tilde w^{\tilde\chi}(c,z),w^{\chi}(c,z))$ est le produit 
$\gamma_1\,\gamma_2\,\gamma_3$ avec
$$\begin{array}{lcl}
\gamma_1 &=& \gamma\big(\tilde\chi^{-1} (w(c,\chi(z))),\chi^{-1} (w(c,\chi(z)))\big)
\\[2pt]
\gamma_2 &=& \gamma(\tilde w(c,\chi(z)),w(c,\chi(z)))
\\[2pt]
\gamma_3 &=& \gamma(\tilde \chi(z),\chi(z))).
\end{array} $$
Soit $L$ un compact de $C$ et $V\in\calv_G$.
Si $W\in\calv_G$ est $G$-invariant et satisfait $W\cdot W\cdot W\subset V$, cela prouve,
en utilisant la formule \eqref{E:inv3}, que si
$(\tilde\chi,\chi)\in\calo^\calg(\alpha(L),W)\cap\calo^\calg(\beta(L),W)$ et 
$(\tilde w,w)\in\calo^\calr(L,W)$, alors 
$(\tilde w^{\tilde\chi},w^{\chi})\in\calo^\calr(L,V)$.  \cqfd

Enfin, l'epace $\calr (\Omega_C,G)$ peut \^etre muni d'une structure
uniforme ${\bf U}_{\calr}$
ayant pour base d'entourages
$$\calo^{\calr}(K,V) :=\{(h,\tilde h)\in\calr (\Omega_C,G)\times\calr (\Omega_C,G)\mid
\tilde h(c)h(c^{-1})\in V ,\ \forall c\in K\}$$
o\`u $V\in \calv_G$ et $K$ est un un compact de $\Omega_C$.

\begin{Proposition} \label{unihom}
Si $G$ est {\it SIN}, l'action 
de $G$ sur $\calr (\Omega_C,G)$
par conjugaison est uniform\'ement continue.
\end{Proposition}

\preu Soit $V\in\calv_G$.
Soit $W\in\calv_G$ tel que $W$ soit  $G$-invariant et $W\cdot W\cdot W\subset V$. 
Soient
$\varphi,\tilde\varphi\in\calr(\Omega_C,G)$ et $g,\tilde g\in G$.
Soit $L$ un compact de $\Omega_C$ et $c\in L$. 
La formule
$$(\tilde g^{-1}\tilde\varphi(c)\tilde g)(g^{-1}\varphi(c) g)^{-1} =
\tilde g^{-1}\tilde\varphi (c)(\tilde g g^{-1})
\tilde\varphi^{-1} (c) (\tilde\varphi (c)\varphi^{-1}(c)) 
\tilde g\,(\tilde g^{-1}   g)$$
montre que $\tilde g\,\tilde g^{-1}\in W$ et
$(\tilde\varphi ,\varphi)\in\calo^\calr (L,W)$ impliquent que \\
$(\tilde g^{-1}\tilde\varphi\tilde g, g^{-1}\varphi g)\in
\calo^\calr (L,V)$.  \cqfd

Nous terminons ce paragraphe en montrant que les 
structures uniformes consid\'er\'ees induisent la CO-topologie. 
Rappelons que, par d\'efinition, une sous-base de la CO-topologie 
sur l'espace fonctionnel ${\rm map\,}(X,Y)$ est form\'ee des ensembles
$CO(K,U):=\{f\in {\rm map\,}(X,Y)\mid f(K)\subset U\}$, o\`u $K$ parcourt
l'ensemble des compacts de $X$ et $U$ celui des ouverts de $Y$.

\begin{Proposition} \label{compou}
Si $G$ est {\it SIN} et $X$ est un espace s\'epar\'e, 
les topologies sur $\calg$, $\calar$ et $\calr(\Omega_C,G)$ 
induites par les structures uniformes  
${\bf U}_\calg$, ${\bf U}_\calar$ et ${\bf U}_\calr$
co\"\i ncident avec la CO-topologie.
\end{Proposition}

Nous aurons besoin d'un analogue du lemme de Lebesgue 
(voir aussi \cite[II.4.3]{Bo})~:

\begin{Lemme}\label{Lebesgue}
Soit $f:K\to E$ une application continue d'un compact dans l'espace total de $\xi$. 
Soit $U$ un ouvert de $E$ avec $f(K)\subset U$. Alors, il existe $V\in\calv_G$ tel
que $f(K)\cdot V\subset U$.
\end{Lemme}

\preu 
Pour tout $z\in K$, il existe un ouvert $S_z$ de $E$ 
et $V_z\in\calv_G$ tels que 
$f(z)\in S_z\subset U$ et $S_z$ est de la forme $\sigma(p(S_z))\times (f(z)\cdot V_z
\cdot V_z)$,
o\`u $\sigma$ est une section locale de $\xi$ au voisinage de $p(f(z))$.
Comme $K$ est compact, on a $f(K)=\bigcup_{z\in K_0}S_z$ pour un sous-ensemble
fini $K_0$ de $K$. Soit $V:=\bigcap_{z\in K_0}V_z\in\calv_G$. 
Alors, pour tout $z\in K$, on a
$f(z)\cdot V\subset U$. En effet, il existe $y\in K_0$ tel que 
$f(z)\in f(y) \dot V_y$ d'o\`u 
$$f(z)\dot V\subset f(y) \cdot V_y\cdot V \subset f(y) 
\cdot V_y\cdot V_y\subset U. \cqfd $$

\vskip .3 truecm\noindent{\sc Preuve de \ref{compou} : }  
L'affirmation pour $\calr(\Omega_C,G)$ est un fait classique pour les applications
dans un espace uniforme \cite[X.3.4, th\'eor\`eme 2]{Bo}.

\vskip .3 truecm\noindent{\it Preuve pour $\calg$ : }
Soit $\chi\in\calg$, $K$ un compact de $E$ et $U$ un ouvert de $E$ 
tels que $\chi(K)\subset U$.
Par le lemme \ref{Lebesgue}, il existe $V\in\calv_G$ tel
que $\chi(K)\cdot V\subset U$.
On a $\chi\in\calo^\calg (p(K),V)(\chi)\subset CO(K,U)$
(observons que $p(K)$ est compact puisque $X$ est s\'epar\'e).

R\'eciproquement, soit $\chi\in\calg$, $L$ un compact de $X$ et $V\in\calv_G$. 
Il faut trouver un ouvert $\Sigma$ pour la CO-topologie tel que 
$\chi\in\Sigma\subset\calo^\calg (L,V)(\chi)$, o\`u
$\calo^\calg (L,V)(\chi):=\{\tilde\chi\in\calg\mid (\chi ,\tilde\chi)\in \calo^\calg (L,V)\}$.
Comme $G$ est \SIN , il existe $W\in\calv_G$ qui est
$G$-invariant avec $W\subset V$. 

Consid\'erons tout d'abord le cas o\`u $L$ est contenu dans un ouvert
$Y$ de $X$ au dessus duquel $\xi$ admet une section $\sigma$. 
L'ensemble $\chi(\sigma(L))\cdot W$ est un ouvert 
de $p^{-1}(L)$. Il existe donc un ouvert $U$ de $E$ tel que
$\chi(\sigma(L))\cdot W=U\cap p^{-1}(L)$. 
Soit $\tilde\chi\in CO(\sigma(L),U)$. Si $z\in p^{-1}(L)$,
il existe $g\in G$ tel que
$z=\sigma(p(z))\cdot g$. La $G$-\'equivariance de $\tilde\chi$ 
donne, pour tout $u\in E$, la formule
\begin{equation}\label{E:Geqchi}
\gamma(\tilde\chi(u\cdot g),\chi(u\cdot g))= 
g^{-1}\,\gamma(\tilde\chi(u),\chi(u))\, g .
\end{equation}
Comme $W$ est $G$-\'equivariant, la formule \eqref{E:Geqchi} appliqu\'ee
\`a $u:=\sigma(p(z))$
entra\^\i ne que $\gamma(\tilde\chi(z),\chi(z))\in W$, 
ce qui prouve que 
$\chi\in CO(\sigma(L),U)\subset\calo^\calg (L,V)(\chi)$.

Dans le cas g\'en\'eral on recouvre $L$ par 
un nombre fini d'ouverts trivialisants,
$L\subset Y_1\cup\cdots Y_n$, au dessus desquels on choisit des sections
$\sigma_i$ de $\xi$. On construit comme ci-dessus les
$\chi(\sigma_i(L_i))\cdot W\subset U_i$ et on aura 
$$\chi\in\bigcap_{i=1}^nCO (\sigma_i(L),U_i)
\subset\calo^\calg (L,W)(\chi)\subset\calo^\calg (L,V)(\chi).  \cqfd $$

\vskip .3 truecm\noindent{\it Preuve pour $\calar$ : }
Soit $w\in\calar$. Soit $L$ un compact de $C\times_X E$ et $U$ un ouvert de $E$ tel que 
$w(L)\subset U$. Par le lemme \ref{Lebesgue}, il existe $V\in\calv_G$ tel que 
$w(L)\cdot V\subset U$ et l'on a 
$w\in \calo^\calr (L,V)(w)\subset CO(L,U)$.

R\'eciproquement, soit $w\in\calar$, $K$ un compact de $C$ et $V\in\calv_G$. 
Comme $G$ est \SIN , il existe $W\in\calv_G$ qui est
$G$-invariant avec $W\subset V$. 
Consid\'erons tout d'abord le cas o\`u $\alpha(K)$ est contenu dans un ouvert
$Y$ de $X$ au dessus duquel $\xi$ admet une section $\sigma$. 
L'ensemble $w(\sigma(\alpha(K)))\cdot W$ est un ouvert 
de $p^{-1}(\beta(K))$. Il existe donc un ouvert $U$ de $E$ tel que
$w(\sigma(\alpha(K)))\cdot W=U\cap p^{-1}(L)$. 
On d\'emontre, comme dans la preuve pour $\calg$ ci-dessus, que
$w\in CO(\sigma(\alpha(K)),U)\subset\calo^\calg (K,V)(w)$.
Le cas g\'en\'eral s'obtient aussi comme dans la preuve pour $\calg$. \cqfd

\section{Groupe de jauge et classifiant de Milnor}  \label{P:gfibuniv}

Les r\'esultats de ce paragraphe seront utilis\'es pour la preuve du
th\'eor\`eme D et les techniques se retrouveront dans la preuve du th\'eor\`eme C.

Soit $\xi:E\hfl{p}{} X$ et $\eta:A\hfl{\bar p}{} B$ deux
$G$-fibr\'e principaux num\'erisables.
Tous les espaces sont point\'es. 
D\'esignons par $\mapers$ l'espace des applications continues point\'ees
$G$-\'equivariantes de $E$ dans $A$, 
muni de la CO-topologie. Le passage au quotient donne
une application $q: \mapers\to\maprs$ o\`u
$\maprs$ est l'espace des applications continues point\'ees
$h:X\to B$ telles que $h^*\eta\approx\xi$, 
muni de la CO-topologie.
Le groupe de jauge $\calg_1$
de $\xi$ agit \`a droite sur $\mapers$ par pr\'e-composition
et $q(f\pcirc\chi)=q(f)$. 

\begin{Theorem}\label{T:gfhinv}
Supposons que $X$ est localement compact,
que $\maprs$ est semi-localement contractile et que $G$ est \SIN .
Alors, l'application 
$q: \mapers\to\maprs$ est un $\calg_1$-fibr\'e principal.
\end{Theorem}

Pour d\'emontrer le th\'eor\`eme \ref{T:gfhinv}, on utilise,
dans l'esprit du \S\ \ref{unif}, 
une sous-base particuli\`ere de la CO-topologie sur $\mapers$. 
Consid\'erons l'ensemble $\calt$ des applications $G$-equivariantes
$\tau : q^{-1}(U_\tau)\to G$ o\`u $U_\tau$ est un ouvert de $B$. Comme
$\eta$ est un $G$-fibr\'e principal, la collection
$\{U_\tau\mid \tau\in\calt\}$ est un recouvrement ouvert de $B$.
Soit $f\in\mapers$.
Soit $K$ un compact de $X$ et $\tau\in\calt$ tel que 
$f(p^{-1}(K))\in q^{-1}(U_\tau)$.
Soit encore $V\in\calv_G$. On d\'efinit
$\bfw (f,K,\tau,V)$ comme l'ensemble des $\tilde f\in\mapers$
telles que ,
pour tout $z\in p^{-1}(K)$ on ait 
$q\pcirc\tilde f(z)\subset U_\tau$ et
$\tau(\tilde f (z)) \in  \tau( f (z))\cdot V$.

\begin{Lemme}\label{T:deutop}
Si $G$ est \SIN , les ensembles $\bfw (f,K,\tau,V)$ 
forment une sous-base de la CO-topologie sur $\mapers$.
\end{Lemme}

\preu Appelons $\bft$ la topologie engendr\'ee par les 
$\bfw (f,K,\tau,V)$ et $\bft_{co}$ la CO-topologie.
Pour $L$ un compact de $E$
et $S$ un ouvert de $A$, notons 
$CO(L,S):=\{h\in\mapers\mid h(L)\subset S\}$. Les
ensembles $CO(L,S)$ forment la sous-base standard de $\bft_{co}$.

Soit $\bfw (f,K,\tau,V)$.  Soit $\sigma:U_\tau\to A$ une section
de $\eta$ restreint \`a $U_\tau$. Comme l'application quotient 
$q(f)\in\maprs$ envoie $K$ dans $U_\tau$, le fibr\'e $\xi$ admet une
section $\hat\sigma$ au dessus de $K$
telle que $f\pcirc \hat \sigma = \sigma\pcirc q(f)$. Soit
$L:=\hat\sigma (K)$ et $S:=\sigma(U)\cdot W$ o\`u $W\in\calv_G$
 est $G$-invariant et contenu dans $V$.
Soit $\tilde f\in CO(L,S)$. Si $z\in p^{-1}(K)$, on a 
$z=z_0\cdot g$ avec $z_0:=\hat\sigma(p(z))\in L$ et
$$\tilde f(z) := \tilde f (z_0)\cdot g \in f(z_0)\cdot W\cdot g
= f(z_0)\cdot g\cdot (g^{-1}\, W\, g) = f(z)\cdot W.$$
Ceci montre que $f\in CO(L,S)\subset \bfw (f,K,\tau,V)$
et donc $\bft \subset \bft_{co}$.

R\'eciproquement, soient $L$ un compact de $E$ et $S$ un ouvert de
$A$ avec $f(L)\subset S$. Supposons tout d'abord que
$S=S(\tau,V):=\sigma(U_\tau)\times V$ pour $V\in\calv_G$
et $\sigma$ la section de $\eta$ restreint \`a $U_\tau$
telle que $\tau\pcirc \sigma (y)=e$, l'\'el\'ement neutre de $G$.
On a donc $\tau\pcirc f (L)\subset V$. Par l'analogue du 
lemme de Lebesgue \cite[II.4.3]{Bo}, il existe $W\in\calv_G$
tel que pour $\tau(f(z))\cdot W\subset V$ pour tout $z\in K$.
Il est alors clair que $f\in W(f,p(L),\tau ,W)\subset CO(L;S)$.

Comme les ouverts $S(\tau,V)$ forment une base 
de la topologie de $A$, on aura, en g\'en\'eral 
$$CO(L,S)=\bigcap _{i=1}^m CO(L_i,S(\tau_i,V_i))\supset 
\bigcap _{i=1}^m  \bfw (f,p(L_i),\tau_i ,W_i) \ni f. $$
On a ainsi prouv\'e $\bft \supset \bft_{co}$. \cqfd

\paragraph{Preuve du th\'eor\`eme \ref{T:gfhinv}~: }

\begin{ccote}\label{PP0} Principe de la d\'emonstration~: \rm 
on va montrer que~:
\begin{enumerate}
\item $q$ est continue, $\calg_1$-\'equivariante et induit une injection 
de $\mapers /\calg_1$ dans $\maprs$.
\item l'action $\mapers\times\calg_1\to\mapers$ est libre et continue.
\item si $f_1,f_2\in\mapers$ satisfont $q(f_1)=q(f_2)$, les points 1) et 2) donnent
un unique $\delta(f_1,f_2)\in \calg_1$ tel que $f_2=f_1\pcirc \delta(f_1,f_2)$.
Ceci d\'efinit une application $\delta : \mapers\times_{\maprs} \mapers\to\calg_1$.
On d\'emontre que $\delta$ est continue.
\item l'application $q$ admet des sections locales continues.
\end{enumerate}

Les points ci-dessus permettent de d\'emontrer que $q$ est un $\calg_1$-fibr\'e principal.
En effet, pour construire des trivialisations locales, on choisit une section continue
$s : T\to \mapers$ de $q$ au dessus d'un ouvert $T$ de $\maprs$. 
La correspondance $(f,\chi)\to s(f)\pcirc \chi$ donne une bijection continue
$\calg_1$-\'equivariante $\tau : T\times\calg_1\to q^{-1}(T)$. Son inverse
$\tau^{-1}(\tilde f) = (q(\tilde f), \delta(s(q(\tilde f)),\tilde f))$ \'etant continue par
3), l'application $\tau$ est un hom\'eo\-mor\-phisme.
\end{ccote}

\begin{ccote}\label{PP1} L'application $q: \mapers\to\maprs$ est continue~: \rm 
Soit $f\in\mapers$. 
Soit $K$ un compact de $X$ et $U$ un ouvert de $B$ avec 
$q(f)\in CO(K,U)$. On peut \'ecrire $K=\bigcup_{i=1}^m K_i$
o\`u $K_i$ est un compact contenu dans le domaine de d\'efinition
d'une section locale $\sigma_i$ de $\xi$.
On a alors 
$$q\big( \bigcap_{i=1}^m CO(\sigma_i(K_i),p^{-1}(U)\big) \subset
\bigcap_{i=1}^m CO(K_i,U) = CO(K,U). $$
\end{ccote}

\begin{ccote}\label{PP2} 
L'action $\mapers\times\calg_1\to\mapers$ est libre et continue~: \rm
Il est banal que cette action est libre.
Soit $(f,\chi)\in\mapers\times\calg_1$. Consid\'erons un 
voisinage de $f\pcirc\chi$ de la forme $\bfw (f\pcirc\chi,K,\tau,V)$.
Soit $W\in\calv_G$ tel que $W\cdot W\subset V$. Comme
dans la preuve de la proposition \ref{unicalg}, on v\'erifie que 
$\tilde f\pcirc \tilde\chi\in \bfw (f\pcirc\chi,K,\tau,V)$
si $\tilde f\in \bfw (f,K,\tau,W)$ et $(\tilde\chi,\chi)\in\calo^\calg(K,W)$.
\end{ccote}

\begin{ccote}\label{PP3} L'application $q$ induit une injection continue
de $\mapers/\calg_1$ dans\\ $\maprs$~: \rm
Il est clair que $q$ est $\calg$-invariante.
Supposons que $q(f_1)=q(f_2)=:f$. On a alors des uniques isomorphismes
$\hat f_i : E\hfl{\approx}{} E(f^*\eta)$, ($i=1,2$).
La composition $\phi:=\hat f_2^{-1}\pcirc \hat f_1$ est un \'el\'ement
de $\calg$ et $f_1=f_2\pcirc \phi$.
(Observons que $q$ est \'evidemment surjective. Nous omettons
ce fait car il est redonn\'e par le point \ref{PP5}. Cette \'economie
sera avantageuse dans la preuve du th\'eor\`eme B).
\end{ccote} 

\begin{ccote}\label{PP4} Continuit\'e de l'application $\delta$ :
\rm Soient $f_1,f_2\in\mapers$ telles que
$f_2=f_1\pcirc \chi$. Il s'agit de montrer que $\chi$ d\'epend contin\^ument
du couple $(f_1,f_2)$. Soient $K$ un compact de $X$ et $V\in\calv_G$.
Soit $W\in\calv_G$ 
tel que $W\cdot W\subset V$.
Supposons tout d'abord qu'il existe $\tau$ telle que 
$K\subset f_1^{-1}(U_\tau)\cap f_2^{-1}(U_\tau)$.
Soient $\tilde f_i\in \bfw (f_i,K,\tau ,W)$ ($i=1,2$)
avec $\tilde f_2 = \tilde f_1\pcirc\tilde\chi$.
Pour $z\in p^{-1}(K)$, on a
$$\tau (\tilde f_2(z)) \in \tau (f_2(z))\cdot W =
\tau (f_1(\chi(z)))\cdot W \in \tau (\tilde f_1(\chi(z)))\cdot W\cdot W
.$$
D'autre part~:
$$\tau (\tilde f_2(z)) = \tau (\tilde f_1(\tilde\chi(z)))=
\tau (\tilde f_1(\chi(z)))\gamma(\tilde\chi(z),\chi(z)).$$
Comme $V=V^{-1}$, on aura $\tilde\chi (z)\in\chi (z) \cdot V$.
Dans le cas g\'en\'eral, on utilise que
$K$ est une r\'eunion finie de compacts $K_i$ tels que 
$K\subset f_1^{-1}(U_{\tau_i})\cap f_2^{-1}(U_{\tau_i})$.
\end{ccote}

\begin{ccote}\label{PP5} Construction de sections locales~: \rm
Nous allons construire une section locale au dessus
de chaque ouvert $T$ de $\maprs$ dont l'inclusion $T\subset\maprs$ est
contractile. On suppose donc qu'il existe une application
$H:I\times T\to \maprs$, o\`u $I=[0,1]$, telle que  
$H(0,f)=f$ et $H(1,f)=f_1$. Comme $X$ est localement compact, $H$ donne 
naissance \`a une application continue $h:I\times T\times X\to B$
\cite[p. 261]{Du} telle que $h(t,f,*)=*$.
D\'esignons par $h_t:T\times X\to A$ l'application
continue $h_t(f,x)=h(t,f,x)$ et par $\cale_t$ l'espace total du
fibr\'e induit sur $T\times X$ par $h_t$~: $\cale_t:=E(h_t^*\eta)$. 
Vu que $h_t(f,*)=*$, l'espace $\cale_t$ est point\'e par $(*,\tilde *)$.

Comme $\eta$ est num\'erisable, le rel\`evement des homotopies donne 
un isomorphisme de $G$-fibr\'es point\'es de
$h_1^*\eta$ sur $h_0^*\eta$. D'autre part, 
on a des isomorphismes de $G$-fibr\'es point\'es
$h_1^*\eta\approx T\times f_1^*\eta\approx T\times \xi$.   
Tout ceci forme un diagramme commutatif
$$\begin{array}{cccccccccccccc}
T\times E & \hfl{\approx}{} & T\times E(h_1^*\eta) &\hfl{\approx}{} &
\cale_1 &\hfl{\approx}{} & \cale_0 &\hfl{}{} & A \\
\downarrow && \downarrow && \downarrow && \downarrow && \downarrow \\
T\times X & \hfl{\rm id}{} & T\times X &\hfl{\rm id}{} &
T\times X &\hfl{\rm id}{} & T\times X &\hfl{h_0}{} & B.
\end{array}$$
Comme $\mapers$ est munie de la CO-topologie,
la ligne sup\'erieure
 d\'eter\-mine une application continue 
$T\to\mapers$ au dessus de l'inclusion
$T\subset \maprs$, c'est-\`a-dire une section locale continue au dessus de $T$. 
\end{ccote}

Par \ref{PP5}, la d\'emonstration du th\'eor\`eme \ref{T:gfhinv} est ainsi termin\'ee. 
De la m\^eme mani\`ere, on d\'emontre le r\'esultat analogue pour
les applications non-point\'ees $q: {\rm Map\,}_G(E,A)\to{\rm Map\,}(X,B)_\xi$~:

\begin{Theorem}\label{T:gfhinvnp}
Supposons que $X$ est localement compact,
que ${\rm Map\,}(X,B)_\xi$ 
est semi-localement contractile et que $G$ est \SIN .
Alors, l'application \\ 
$q: {\rm Map\,}_G(E,A)\to {\rm Map\,}(X,B)_\xi$ est un $\calg$-fibr\'e principal.
 \end{Theorem}

Le cas particulier o\`u $\eta$ est le fibr\'e de Milnor $EG\to BG$ est
int\'eressant \`a cause de la proposition suivante, utilis\'ee pour d\'emontrer le th\'eor\`eme D~:

\begin{Proposition}\label{T:mapegcontr}
Soit $\xi$ un $G$-fibr\'e principal num\'erisable sur $X$.
Alors, l'espace
${\rm Map\,}_G(E,EG)$  est contractile. Si, de plus, 
$X$ est localement compact et que
l'inclusion $\{*\}\subset X$ soit une cofibration,
${\rm Map\,}_G^{\scriptscriptstyle\bullet}(E,EG)$
est faiblement contractile (i.e. ses groupes d'homotopie sont ceux d'un point).
\end{Proposition}

\preu La preuve habituelle que deux applications $G$-equivariantes
$f_0,f_1 : E\to EG$ sont toujours homotopes fournit une homotopie
{\it canonique} entre $f_0$ et $f_1$, qui d\'epend contin\^ument de $(f_0,f_1)$
(voir \cite[prop. 12.3]{Hu}).
L'espace ${\rm Map\,}_G(E,EG)$ est donc contractile (convexe), de m\^eme que 
$EG = {\rm Map\,}_G(G,EG)$.
Pour montrer l'assertion sur 
${\rm Map\,}_G^{\scriptscriptstyle\bullet}(E,EG)$, il suffit de montrer que 
la suite 
$${\rm Map\,}_G^{\scriptscriptstyle\bullet}(E,EG)\hookrightarrow
{\rm Map\,}_G(E,EG)\hfl{\rm ev}{} EG$$
est, avec nos hypoth\`eses sur $X$, 
une fibration de Hurewicz, o\`u ${\rm ev}$ est l'\'eva\-lua\-tion sur le point base
$\tilde *$.

Soit $f:A\to {\rm Map\,}_G(E,EG)$ une application continue et
$F_{\tilde *}:  [0,1]\times A\times \{\tilde *\}\to EG$ une homotopie de ${\rm ev}\pcirc f$. 
Par passage au quotient, on obtient
$\bar f:A\to {\rm Map\,}(X,BG)$ et $\bar F: [0,1]\times A\times \{*\}\to BG$. Comme $X$ est localement
compact, ces applications induisent des applications continues
$\bar f\,\check{} : A\times X\to BG$ 
et $\bar F_*\check{} : [0,1]\times A\times \{*\}\to BG$.
Etant \'equivariante, l'application induite $f\,\check{} : A\times E\to EG$
est donc aussi continue (bien que $E$ ne soit pas suppos\'e localement compact).

Comme $*\subset X$ est une cofibration, il existe une r\'etraction de $[0,1]\times X$
sur $\{0\}\times X \cup [0,1]\times \{*\}$ qui permet d'\'etendre
$\bar F_*\check{}$ en 
$\bar F\check{} : [0,1]\times A\times X\to BG$. Ceci prouve que 
$ {\rm Map\,}^{\scriptscriptstyle\bullet}(X,BG)\to  
{\rm Map\,}(X,BG)\to BG$ est une fibration de Hurewicz. 
Par rel\^evement des homotopies, l'espace total du
fibr\'e induit sur $[0,1]\times A\times X$ par $\bar F\check{}$ est $G$-hom\'eomorphe
\`a $[0,1]\times A\times E $, ce qui produit une homotopie
$F\tilde {} : [0,1]\times A\to  {\rm Map\,}_G(E,EG)$ partant $f$ et au dessus 
de $F$.  \cqfd

\begin{ccote} \bf Remarque : \rm
On d\'eduit de \ref{T:gfhinv}, \ref{T:gfhinvnp} et \ref{T:mapegcontr} que 
$\pi_{i}({\rm Map\,}^\bullet(X,BG)_\xi)\approx \pi_i(B\calg_1) $ et 
$\pi_{i}({\rm Map\,}(X,BG)_\xi)\approx \pi_i(B\calg) $. 
Nous ne savons pas, en g\'en\'eral, si ces isomorphismes sont induits
par une application ${\rm Map\,}(X,BG)_\xi\to B\calg$. Observons que des
r\'esultats analogues ont \'et\'e obtenus dans d'autres contextes
(\cite{DDK}, \cite[Prop. 5.1.4]{DK}).
\end{ccote}

\section{Preuve des th\'eor\`emes B, C et D}  \label{P:pthmbc}

\paragraph{Pr\'eparatifs : }
\begin{ccote} Continuit\'e des foncteurs de Milnor. \rm 
Nous aurons besoin de savoir que les foncteurs $E$ et $B$ de Milnor sont continus, 
ce qui ne semble pas figurer dans la litt\'erature~:
\end{ccote}
\begin{Lemme}\label{L:contB}
Soient $\calr (F,G)$ l'espace des homomorphismes continus
entre les groupes topologiques $F$ et $G$. Supposons que $F$ est s\'epar\'e.
Alors les applications
$$ E :  \calr (F,G)\to {\rm Map\,}_F^{\scriptscriptstyle\bullet}(EF,EG)
 \qquad ( \phi\ \mapsto \ E\phi )$$ et
$$ B :  \calr (F,G) \to {\rm Map\,}^{\scriptscriptstyle\bullet}(BF,BG) 
\qquad ( \phi\ \mapsto \ B\phi )$$
sont continues (tous les espaces \'etant munis de la CO-topologie).
\end{Lemme}

\preu 
Comme l'application 
${\rm Map\,}_F^{\scriptscriptstyle\bullet}(EF,EG)\to
{\rm Map\,}^{\scriptscriptstyle\bullet}(BF,BG)$
est continue (se d\'emontre comme \ref{PP1}), 
il suffit de donner une preuve pour l'application $E$. 

Rappelons que si $P$ est un groupe toplogique, les \'el\'ements de $EP$ 
sont repr\'esent\'es par des suites $(t_jp_j)$, o\`u $j\in\bbn$,
$t_j\in [0,1]$ et $p_j\in P$. 
On d\'esigne par $\tau_i : EP\to [0,1]$ l'application
$\tau_i((t_jp_j)):=t_i$ et par $\gamma_i: \tau_i^{-1}(]0,1])\to P$ 
l'application  $\gamma_i((t_jp_j)):=p_i$
(on utilise les m\^emes notations $\tau_i$ et $\gamma_i$ pour tout groupe $P$). 
L'espace $EP$ est muni de la topologie 
la plus grossi\`ere telle que les applications $\tau_i$ et $\gamma_i$ soient continues. 
Une sous-base $\cals$ de cette topologie est donc constitu\'ee par les ouverts du type
\begin{enumerate} 
\item\label{ty1} $\tau_i^{-1}(J)$ o\`u $J$ est un ouvert de $[0,1]$.
\item\label{ty2} $\gamma_i^{-1}(V)$ o\`u $V$ est un ouvert de $P$.
\end{enumerate} 

Soit $\phi_0\in\calr (F,G)$,  
$K$ un compact de $EF$ et $U$ un ouvert de $EG$ tel que $\phi_0(K)\subset U$.
Soit $CO(K,U):=\{\alpha\in {\rm Map\,}_F^{\scriptscriptstyle\bullet}(EF,EG) \mid
\alpha(K)\subset U\}$. Il faut montrer que $E^{-1}(CO(K,U))$ est un voisinage de
$\phi_0$ dans $\calr (F,G)$. Il suffit de le faire pour $U\in \cals$
(voir \cite[X.3.4, remarque 2]{Bo}. 

Si $U$ est du type \ref{ty1} ci-dessus, cela ne pose pas de probl\`eme. En effet,
$\tau_i\pcirc E = \tau_i$, d'o\`u $E^{-1}(CO(K,U))=\calr (F,G)$.
Si $U=\gamma_i^{-1}(V)$ avec $V$ un ouvert de $G$, on pose $K_i:=\gamma_i(K)$
($\gamma_i$ est d\'efinie sur $K$ puisque $\tau_i\pcirc E = \tau_i$).
L'espace $K_i$ est compact puisque $F$ est s\'epar\'e.
Comme $\gamma_i\pcirc E\phi = \phi\pcirc\gamma_i$, on a 
$\phi_0\in CO(K_i,V) \subset E^{-1}(CO(K,U))$.  \cqfd

\begin{ccote}\label{applicn}  L'applcation $n :\calrr\to\mapr$. \rm
En choisissant une partition de l'unit\'e $\mu_\rho$ subordonn\'ee au recouvrement
$U_\rho$ ($\rho\in {\rm Cont\,}_C(X)$), on d\'etermine, gr\^ace aux
trivialisations $\hat\psi_\rho$ de \eqref{deftriv2},
une application classifiante $\nu_C:X\to B\Omega_C$ 
(voir \cite[ch. 4, prop. 12.1 et th. 12.2]{Hu}).
Il est possible de faire ce choix de mani\`ere que
$\nu_C$ soit point\'ee (le point base de $EF$, 
pour un groupe topologique $F$,
est toujours la suite $(1e,0,0,\dots)$, o\`u $e$ \'el\'ement neutre de $F$; 
le point base de $BF$ est l'image de celui de $EF$).
Pour cela, il faut tout d'abord avoir une partition
de l'unit\'e d\'enombrable et trivialisante $\mu_i$ ($i=1,2,\dots $) 
telle que $\mu_1(*)=1$. 
Or, une telle partition existe car~:

a) il existe une partition de l'unit\'e $\hat \mu_\rho$ et $\hat\rho\in
{\rm Cont\,}_C(X)$ tels que 
$\hat \mu_{\hat\rho(x)}=1$ au voisinage de $*$. Pour voir cela,
on choisit $\hat\rho$ tel que $\mu_{\hat\rho}(*)\not = 0$. On consid\`ere
une fonction $\delta: [0,1]\to [0,1]$ telle que 
$\delta(t) = 0$ si $t\geq \mu_{\hat\rho}(*)/2$ et
$\delta(t) > 0$ si $t < \mu_{\hat\rho}(*)/2$. On pose
$$\mu'_\rho(x) := \left\{
\begin{array}{lll} \mu_{\hat\rho}(x) & \hbox{si}\
\rho = \hat\rho\\
\delta(\mu_{\hat\rho}(x))\,\mu_\rho(x) & \hbox{sinon.}\end{array}\right.$$
Avec ces d\'efinitions, 
$$\hat\mu_\rho (x):=\frac{\mu'_\rho(x)}{\sum_{\sigma}\mu'_\sigma(x)}$$
est une partition de l'unit\'e avec 
$\hat \mu_{\hat\rho(x)}=1$ au voisinage de $*$.

\vskip .2 truecm
b) le proc\'ed\'e de \cite[ch. 4, prop. 12.1]{Hu} fournit, \`a partir
de la partition $\hat\mu_\rho$ une partition de l'unit\'e
$\mu_i$ ($i=1,2,\dots $) telle que $\mu_1(x)=1$ au voisinage de $*$.

On peut supposer que 
$\hat\rho(*)=i_{\hat\rho(*)}$. Sinon, $\hat\rho(*)=b\in\Omega_C$
et l'on remplace $\hat\rho$ par $b^{-1}\,\hat\rho$.
Dans ces conditions, l'application $\nu_C$ obtenue 
par \cite[ch. 4, prop. 12.1]{Hu} est point\'ee.

Ayant fix\'e $\nu_C\in {\rm map\,}(X,B\Omega_C)$ comme ci-dessus,
on d\'efinit $n: \calrr\to\mapr$ par
$n(\phi):=B\phi\pcirc \nu_C$. Comme $C$ est s\'epar\'e, 
l'application $\phi\mapsto B\phi$ est 
continue (lemme \ref{L:contB}), d'o\`u la continuit\'e de  $n$. 
\end{ccote}

\paragraph{Preuve du th\'eor\`eme B : }
Par le th\'eor\`eme d'existence, on a $\calar=\emptyset$ si et seulement 
$\calr (\Omega_C,G)_\xi = \emptyset$. 
Dans ce cas, le th\'eor\`eme B est banal.
On suppose donc que $\calar\not =\emptyset$.
Le th\'eor\`eme A implique que 
l'image de $h$ est dans $\calr (\Omega_C,G)_\xi$.
Le principe de la preuve du th\'eor\`eme B est alors le m\^eme que celui du th\'eor\`eme
\ref{T:gfhinv} (voir \ref{PP0}).

\begin{ccote}  \label{hconti} $h$ est continue. \ \rm 
Soit $w\in\calar$. Soit $K$ un compact de $\Omega_C$, $U$ un ouvert de $G$ tels que
$h^w(K)\subset U$. En utilisant une trivialisation locale de $\xi$ au voisinage
de $*$, on peut trouver un ouvert $\tilde U$ de $E$ tel que 
$\tilde U \cap E_* = \tilde * \cdot U$. 
Alors $w\in CO(K\times\{\tilde *\},\tilde U)$ qui est un ouvert de $\calar$ et 
\\ $h(CO(K\times\{\tilde *\},\tilde U)\subset CO(K,U)$.
\end{ccote}

\begin{ccote}  \label{L:g1lib}
l'action $\calar\times\calg_1\to\calar$ est continue et libre. \ \rm
La continuit\'e, utilisant le fait que $G$ est \SIN , a \'et\'e d\'emontr\'ee dans la proposition \ref{unicala}.
Soit $w\in\calar$ et $\chi\in\calg$. Si $\chi\not ={\rm id}_E$, il
existe $z\in E$ tel que $\chi(z)\not = z$. 
Soit $\theta\in C_{\scriptscriptstyle p(z)}^*$.
Comme $\chi$ restreinte
\`a $p^{-1}(*)$ est l'identit\'e, on aura 
$w^\chi (\theta,z)\not = w(\theta,z)$ 
ce qui montre l'assertion \ref{L:g1lib}. 
\end{ccote}

\begin{ccote} $\calg$-invariance de $h$ : \  \rm
Soit $\chi\in\calg$. Posons $\chi(\tilde *) = \tilde *\, g$,
autrement dit~: $g:=\res(\chi)$. 
Pour $c\in\Omega_C$, on a
$$\begin{array}{lcl}
\tilde *\, h^{(w^\chi)}(c) & = & w^\chi(c,\tilde *) = 
\chi^{-1}(w(c,\tilde * \cdot g)) = \chi^{-1}(\tilde * 
\cdot h^w(c))\cdot g =\\[2pt] & = &
\chi^{-1}(\tilde * )\cdot (h^w(c)\, g) = \tilde *\cdot (g^{-1}\, h^w(c)\, g)
\end{array}$$
d'o\`u 
\begin{equation}\label{conju34}
h^{(w^\chi)}(c)= \res (\chi)^{-1}\, h^w(c)\, \res(\chi)
\end{equation}
et $h$ induit des
applications
$$\bar h_1: \calar/\calg_1\to\calr (\Omega_C,G)_\xi \quad \hbox{ et } \quad
\bar h: \calar/\calg\to\calr (\Omega_C,G)_\xi/G.$$
\end{ccote}

\begin{ccote}\label{injbah1}  {\it Injectivit\'e de $\bar h_1$ : } \rm
Soient $w,\tilde w\in\calar$ telles que $h^w=h^{\tilde w}$.  
Soit $z\in E$. Choisissons $\rho\in {\rm Cont\,}_C(X)$ tel que 
$p(z)\in U_\rho$ et notons $\theta:=\rho(p(z))\in C_{p(z)}^*$. 
On d\'efinit 
\begin{equation}\label{trgau}
\chi (z):= 
\tilde w (\theta^{-1},w(\theta,z)).
\end{equation}
Nous allons montrer que l'\'egalit\'e $h^w=h^{\tilde w}$entra\^\i ne
que $\chi(z)$ 
ne d\'epend pas du choix 
de  $\rho$. 
Soit $\bar \rho$ une autre $C$-contraction, donnant 
$\bar\theta$ et $\bar \chi (z)$. Rappelons que la fibre $p^{-1}(*)$
est identifi\'ee \`a $G$ par $g\mapsto \tilde *\cdot g$. Via cette identification, $G$ agit \`a gauche
sur $p^{-1}(*)$ et, si
$c\in\Omega_C$ et $y\in p^{-1}(*)$, on a $w(c,y)= h^w(c)\cdot y$. Avec
ces conventions, on a
$$\begin{array}{lcl}
\bar\chi(z) &=& \tilde w (\bar\theta^{-1},w(\bar\theta,z)) =
\tilde w (\bar\theta^{-1},w(\bar\theta\theta^{-1}\theta,z)) = \\[2pt]&=&
\tilde w (\bar\theta^{-1},w(\bar\theta\theta^{-1}, w(\theta,z))) =
\tilde w (\bar\theta^{-1},h^w(\bar\theta\theta^{-1})\cdot w(\theta,z)) =
\\[2pt]&=&
\tilde w (\theta^{-1}\theta\bar\theta^{-1},h^w(\bar\theta\theta^{-1})\cdot w(\theta,z)) =
\\[2pt]&=&
\tilde w (\theta^{-1},
\underbrace{h^{\tilde w}(\theta\bar\theta^{-1})
h^w(\bar\theta\theta^{-1})}_{= 1}\cdot w(\theta,z)) =
\tilde w(\theta^{-1},w(\theta,z))
= \chi(z).
\end{array}$$
On a ainsi d\'efini une application $\chi : E\to E$ qui, par la formule
\eqref{trgau} est continue. Son inverse s'obtient en \'echangeant  
$w$ et $\tilde w$. Les formules
$\chi(z\cdot g) = \chi(z)\cdot g$, pour $g\in G$ et 
$p(\chi(z))=p(z)$ sont banales. 
De plus, on a $\chi(\tilde *)=\tilde *$,
d'o\`u $\chi\in\calg_1$.

Voyons maintenant que $\tilde w^\chi = w$. Soit $z\in E$ et  
$c\in C_{p(z)}$.
Observons que, dans la
la formule \eqref{trgau}, on n'utilise la $C$-contraction 
$\rho$ que pour garantir la continuit\'e. La d\'efinition de
$\chi(z)$ ne n\'ecessite que l'\'el\'ement $\theta\in C_{p(z)}^*$ et 
$\chi(z)$ ne d\'epend pas de $\theta$. En choisissant $\theta$ pour la
d\'efinition de $\chi(z)$ et $\theta c^{-1}$ pour celle de 
$\chi^{-1}(\tilde w(c,\chi(z))$, on aura 
$$\begin{array}{lcl}
\tilde w^\chi(c,z) &=& \chi^{-1}(\tilde w (c,\chi(z))) =
\chi^{-1}(\tilde w (c,\tilde w(\theta^{-1},w(\theta,z))) = \\[2pt]&=&
\chi^{-1}(\tilde w (c\theta^{-1},w(\theta,z)) =
w(c\theta^{-1},\tilde w(\theta c^{-1},\tilde w (c\theta^{-
1},w(\theta,z))= \\[2pt]&=&
w(c\theta^{-1},w(\theta,z)) = w(c,z). 
\end{array}$$
\end{ccote}

\begin{ccote} 
L'application $\delta : \calar\times_{\calrr}\calar\to\calg$ 
est uniform\'ement continue~: \rm 
Soit $V\in\calv_G$ et $K$ un compact de $X$.  
Supposons d'abord que $K\subset U_\rho$ pour 
une $C$-contraction $\rho\in {\rm Cont\,}_C(X)$.
Comme $C$ est s\'epar\'e, les sous-espaces
$L$ et $L^{-1}$ de $C$ d\'efinis par 
$$L:= \{\rho(x)\mid x\in K\}\subset C^* \qquad \hbox{et} \qquad
L^{-1}:= \{\rho(x)^{-1}\mid x\in K\}\subset C_*.$$
sont compacts. Soit $W\in\calv_G$ tel que $W\cdot W\in V$. 
Soient $\tilde w_1,\tilde w_2\in \calar$.
Il suit de \ref{injbah1} que $\delta$ satisfait, pour 
$z\in p^{-1}(K)$, \`a l'\'equation
\begin{equation}\label{E:pfb111}
\delta(\tilde w_1,\tilde w_2)(z) 
= \tilde w_1(\theta^{-1},\tilde w_2(\theta,z))
\end{equation}
avec $\theta:=\rho(p(z))\in C_{p(z)}^*$.
Il s'en suit que si 
$(\tilde w_1, w_1)\in \calo^\calr (L^{-1},W)$ et 
$(\tilde w_2, w_2)\in \calo^\calr (L,W)$, alors
$(\delta(\tilde w_1,\tilde w_2),
\delta(w_1,w_2))\in \calo^{\calg_1}(K,V)$.

Dans le cas g\'en\'eral, on utilise que $K:=\bigcup_{\rho\in P} K_\rho$
o\`u $P$ est un ensemble fini dans ${\rm Cont\,}_C(X)$ et
$K_\rho$ est un compact de $U_\rho$. On d\'efinit, comme ci-dessus,
$L_\rho:=\rho(K_\rho)$ et $L_\rho^{-1}:=\rho(K_\rho)^{-1}$ et on aura
$$\left.\begin{array}{l}\displaystyle
(\tilde w_1, w_1)\in \bigcap_{\rho\in P}\calo^\calr (L_\rho^{-1},W)\\
\hbox{et}\\[1 pt] \displaystyle
(\tilde w_2, w_2)\in \bigcap_{\rho\in P}\calo^\calr (L_\rho,W)
\end{array}\right\} \ \Rightarrow \ 
(\delta(\tilde w_1,\tilde w_2),
\delta(w_1,w_2))\in \calo^{\calg_1}(K,V).$$
\end{ccote}

\begin{ccote}\label{trivlocthb}  Sections locales : \rm
Soit $T$ un ouvert de $\mapr$ tel que l'inclusion $T\subset\mapr$ soit
contractile. Soit  $S:=n^{-1}(T)$, o\`u $n:\calrr\to\mapr$ est
l'application d\'efinie en \ref{applicn}.
Consid\'erons l'application compos\'ee
\begin{equation}\label{defN1}
N: \, S\times X   \hfl{n\times {\rm id}}{} T\times X 
\hfl{{\rm ev}_X}{}  BG
\end{equation}
Comme $X$ est localement compact, l'\'evaluation ${\rm ev}_X$ est continue 
et $N$ est continue. D\'esignons par $\cale:= N^*EG$, l'espace total du 
$G$-fibr\'e principal sur $S\times X$ induit par $N$. 

Observons que $N$ est aussi la composition 
\begin{equation}\label{defN2}
N : S\times X   \hfl{{\rm id_S}\times \nu_C}{} S\times 
B\Omega_C \hfl{{\rm ev}_\Omega}{} BG
\end{equation}
de ${\rm id_S}\times \nu_C$ avec l'application d'\'evaluation
(peut-\^etre non-continue) ${\rm ev}_\Omega$.
Comme le $\Omega_C$-fibr\'e induit sur $X$ par  
$\nu_C$ est $\xi_C : C_*\hfl{\beta}{} X$, 
on obtient, en utilisant \eqref{defN2}, que l'espace 
$\cale$ est le quotient de $S\times C_*\times G$ par la relation 
d'\'equivalence $(\phi ,ub,g)\sim (\phi ,u, \phi(b)g)$, o\`u
$(\phi ,u,g)\in S\times C_*\times G$ et $b\in\Omega_C$.
Consid\'erons $\cale$ comme un espace au dessus de $X$ par l'application
$(\phi ,u,g)\mapsto \beta (u)$ et 
d\'efinissons l'application continue 
$v: C\times _X \cale \to \cale$ par $v(c,(\phi ,u,g)):=(\phi ,cu,g)$.

Comme l'application $S\to\mapr$ est homotope \`a une application constante, 
on montre, comme dans \ref{PP5},
qu'il existe
un hom\'eomorphisme $G$-equivariant $S\times E \hfl{\approx}{}\cale$
au dessus de ${\rm id}_{S\times X}$. Via cet hom\'eomorphisme
et en composant avec la projection $S\times E\to E$, l'application
$v$ donne une application continue 
$\hat v: C\times_X (S\times E)\to E$. A son tour, $\hat v$ d\'etermine
une application continue $w : S\to {\rm map\,}(C\times_X E, E)$.
On v\'erifie facilement que l'image de $w$ 
est dans $\calar$. La pr\'eservation des points base par les divers isomorphismes
utilis\'es \cb
entra\^\i ne que $w$ est une section locale de $h$ au dessus de $S$. 
\end{ccote}

Ayant \'etabli les points \ref{hconti} \`a \ref{trivlocthb}, la d\'emonstration du
th\'eor\`eme B se termine comme expliqu\'e dans \ref{PP0}. \cqfd

\paragraph{Preuve du th\'eor\`eme C : }
Soient $w,\tilde w\in\calar$ et $g\in G$ tels que 
 $h^w=g^{-1}h^{\tilde w}g$. Comme $G$ est connexe par arc, il existe, par
le lemme \ref{numer}, un \'el\'ement
$\chi\in\calg$ tel que $\chi(\tilde *)=\tilde * \cdot g$.
L'\'equation \eqref{conju34} montre qu'alors
$h^{w^\chi}= h^{\tilde w}$. Comme $\bar h_1$ est injective
par \ref{injbah1}, les repr\'esentations ${w^\chi}$ et $\tilde w$ 
sont dans la m\^eme classe modulo $\calg_1$. Cela prouve 
l'injectivit\'e de $\bar h$. 
Le diagramme commutatif
\begin{equation}\label{hh1}
\begin{array}{cccccc}
\calar/\calg_1 & \hfl{\bar h_1}{} & \calr (\Omega_C,G)_\xi \\
\downarrow && \downarrow \\
\calar/\calg & \hfl{\bar h}{} & \calr (\Omega_C,G)_\xi/\hbox{\small 
G}
\end{array}\end{equation}
et le fait que $\bar h_1$ soit un hom\'eomorphisme (vu le th\'eor\`eme B)
font que $\bar h$ est un hom\'eomorphisme.

\paragraph{Preuve du th\'eor\`eme D :}\label{P:pthmd}
\begin{Lemme}\label{appator}
L'application $n$ est couverte par un morphisme de
$\calg_1$-fibr\'es principaux
\begin{equation}\label{appatordiag}
\begin{array}{cccccc}
\calar\ \ & \hfl{\hat n}{} & \maper \\
\downarrow\! h && \downarrow\\
\calrr  & \hfl{n}{} & \mapr
.\end{array}\end{equation}
\end{Lemme}

Le th\'eor\`eme D d\'ecoulera du lemme \ref{appator} puisque,
i l'inclusion $\{*\}\subset X$ est une cofibration, l'espace
$\maper$ est faiblement contractile (Proposition \ref{T:mapegcontr}).
En fait, l'application $n$ est classifiante pour le $\calg_1$-fibr\'e
principal $h : \calar\to\calrr$.

\vskip .2 truecm
{\sc Preuve du lemme \ref{appator} : } 
Le proc\'ed\'e pour fixer l'application $\nu_c$ vu en 
\ref{applicn} produit en fait un diagramme commutatif
$$\begin{array}{cccc}
C_* & \hfl{\hat \nu_C}{} & E\Omega_C\\
\downarrow\!\! \beta &&\downarrow\\
X & \hfl{\nu_C}{} & E\Omega_C\end{array}$$
o\`u 
$\hat \nu_C\in {\rm map\,}_G^{\scriptstyle\bullet}(C_*,E\Omega_C)$.
L'application $\hat n\in 
{\rm map\,}_\calg(\calar,\maper)$ est d\'efinie
par $\hat n(w):=Eh^w\pcirc \hat\nu_C$. Le diagramme 
\ref{appatordiag} est bien commutatif.
Observons que $\hat n$ est obtenue par le proc\'ed\'e de 
\cite[ch. 4, prop. 12.1]{Hu}
\`a l'aide de la partition de l'unit\'e
$\hat \mu_\rho$ construite en \ref{applicn} et des trivialisations
$\psi^w_\rho$ du lemme \ref{rectrican}.

\section{Exemples et applications} \label{P:exp}

\subsection{Le groupo\"\i de associ\'e \`a un fibr\'e principal} \label{SP:groassfp}

Soit $\xi : E\hfl{p}{} X$ un $G$-fibr\'e principal. Soit $EE:=EE(\xi):=(E\times E)/G$, le quotient
de $E\times E$ par l'action diagonale de $G$. D\'enotons par 
$\llangle{a}{b}$ l'orbite de $(a,b)$ dans $EE$. On fait de $EE$ un $X$-\gro\
en posant $\alpha (\llangle{a_2}{a_1}) := p(a_1)$, $\beta (\llangle{a_2}{a_1}) := p(a_2)$,
$i_x:=\llangle{a}{a}$ avec $p(a)=x$; la composition vient de la formule
$\llangle{a_3}{a_2}\llangle{a_2}{a_1}=\llangle{a_3}{a_1}$, ou, plus g\'en\'eralement,
$$\llangle{a_3}{a'_2}\llangle{a_2}{a_1}=\llangle{a_3\cdot \gamma(a'_2,a_2)}{a_1}$$
o\`u $\gamma$ est l'application d\'efinie en \eqref{defgamma}. L'inverse est \'evidemment donn\'e
par $\llangle{a}{b}^{-1}=\llangle{b}{a}$. Le $X$-\gro\ ainsi obtenu s'appelle
le {\it $X$-\gro\ associ\'e \`a $\xi$} \cite[p. 5]{Ma}.

Le groupe structural $G$ de $\xi$ est isomorphe \`a $\Omega_{EE}$ par 
$g\mapsto \llangle{\tilde * \cdot g}{\tilde *}$. Observons que $EE$ est localement
trivial num\'erisable si $\xi$ est num\'eridable. En effet, 
$x\mapsto \llangle{\tilde *}{\sigma(x)}$ est une $EE$-contraction lorsque $\sigma$ est
une section locale de $\xi$. 

Le $X$-\gro\ $EE$ a une repr\'esentation tautologique $v$ sur $\xi$ d\'etermin\'ee par 
l'application continue $v : (E\times E)\times E \to E$ d\'efinie par
$v((a,b),z):=a\cdot\gamma(b,z)$; son holonomie est ${\rm id\,}_G$. 
Pour un $X$-\gro\ $C$, d\'esignons par ${\rm Mor\,}(C,EE)$ l'ensemble des
morphismes continus de $C$ dans $EE$ (au dessus de ${\rm id\,}_X$).

\begin{Proposition}\label{reptoto}
La composition avec la repr\'esentation tautologique donne une bijection
de ${\rm Mor\,}(C,EE(\xi))$ sur $\calar$.
\end{Proposition} 

\preu La bijection inverse $w\mapsto \Psi_w$ est donn\'ee par 
$\Psi_w(c):=\llangle{w(c,z)}{z}$ ($z\in E_{\alpha(c)}$). La seule chose non-triviale 
\`a v\'erifier est que $\Psi_w$ est continu. 

Soit $c\in C$ et $T\ni \Psi_w(c)$ un ouvert de $EE$. Soit $z\in E_{\alpha(c)}$. 
Il existe 
\begin{itemize}
\item $V$ un ouvert de $E$ contenant $w(c,z)$,
\item  $U$ un ouvert de $X$ contenant $x:=\alpha(c)=p(z)$ et $\sigma :U\to E$ une section
continue locale de $\xi$ avec $\sigma(x)=z$ et 
\item $W\in \calv_G$
\end{itemize}
tels que $\pi (V\times (\sigma (U)\cdot W))\subset T$, o\`u $\pi$ d\'esigne la projection
de $E\times E$ sur $EE$. Par continuit\'e de $w$, il existe 
$A$ un ouvert de $C$ contenant $c$ et $R$ un ouvert de $E$ contenant $z$ tels que
$w(A\times R)\subset V$. Soient $U'$ un ouvert de $X$ et $W'\in\calv_G$ tels que
$x\in U'\subset U$, $W'\subset W$ et $\sigma(u')\cdot W'\subset R$.
Soit $A':= A\cap \alpha^{-1}(U')$. On a $c\in A'$ et $\Psi_w(A')\subset 
\pi(V\times (\sigma(U')\cdot W')\subset T$, ce qui prouve la continuit\'e de $\Psi_w$ en $c$.
\cqfd

\subsection{Pr\'egroupo\"\i des} \label{SP:pregro}

Comme nous le verrons plus loin, une fa\c con commode de d\'efinir un \gro\ 
topologique est de partir d'une cat\'erorie
topologique avec anti-involution (pr\'e\gro ).
Soit $X$ un espace topologique muni d'un point base $*\in X$. 
Un $X$-{\it pr\'e\gro\ } est un espace topologique
$\tilde C$ muni de deux applications continues $\alpha,\beta : \tilde C\to X$, 
d'une {\it composition} partiellement d\'efinie et d'une application
$x\mapsto i_x$ de $X$ dans $\tilde C$ qui satisfont aux propri\'et\'es a), b) et c)
de la d\'efinition du \S\ 1 d'un $X$-\gro .
En revanche, la condition d) consiste seulement en 

d') $\tilde C$ est muni d'une anti-involution continue $c\mapsto \bar c$, 
envoyant $\tilde C^y_x$ sur $\tilde C^x_y$.

Une {\it repr\'esentation} d'un $X$-pr\'e\gro\ $\tilde C$ 
 sur un $G$-espace $\xi$ est 
une application continue $w: \tilde C\times_X E\to E$ 
(o\`u $\tilde C$ est vu au dessus de $X$ via $\alpha$) telle que, pour tout 
$c,d\in \tilde C$, $z\in E$ et $g\in G$, on ait
\begin{enumerate}
\item $p(w(c,z)) = \beta(c)$.
\item $w(c\, d,z)=w(c,w(d,z))$.
\item $w(\bar c,(w(c,z))=z$.
\item $w(c,z\cdot g) = w(c,z)\cdot g$.
\end{enumerate}

Si $X$ est s\'epar\'e, un $X$-pr\'e\gro\ d\'etermine un $X$-\gro\ s\'epar\'e, avec la propri\'et\'e suivante :

\begin{Proposition}\label{predefgro}
Soit $\tilde C$ un $X$-pr\'e\gro\ localement trivial num\'erisable avec $X$ s\'epar\'e. 
Alors, il esiste un unique $X$-\gro\ s\'epar\'e localement trivial num\'erisable
$C$ avec un morphisme continu 
surjectif $\tilde C\to C$ 
satisfaisant \`a la condition suivante~: toute repr\'esentation de 
$\tilde C$ sur un $G$-fibr\'e principal
$\xi$ au dessus de $X$, avec $G$ s\'epar\'e, est induite par une unique 
repr\'esentation de $C$.
\end{Proposition}

La d\'emonstration de \ref{predefgro} utilise le lemme suivant~:

\begin{Lemme}\label{grosep}
Soit $C$ un $X$-\gro\ localement trivial. 
Alors, $C$ est s\'epar\'e si et seulement si $X$ et $\Omega_C$ le sont.
\end{Lemme}

\preu Supposons que $X$ et $\Omega_C$ soient s\'epar\'es (l'autre sens est banal).
Comme $C$ est localement trivial, $\beta : C_*\to X$ est un $\Omega_C$-fibr\'e principal,
par le th\'eorme A. L'espace $C_*$ est donc s\'epar\'e. 
Il en est \'evidemment de m\^eme pour $C^*$.
On utilise alors
que $C= C_*\times_{\Omega_C}C^*$ pour \'etablir que $C$ est s\'epar\'e. \cqfd

\vskip .2 truecm\noindent
{\sc Preuve de la proposition \ref{predefgro} : } Soit $\hat C$ l'ensemble quotient de 
$\tilde C$ par la relation d'\'equivalence engendr\'ee par 
$c\bar u u d \sim cd$, pour tout 
$c,d,u\in\tilde C$ avec $\beta (d) = \alpha(u) = \alpha (c)$.
Il est clair que la structure (alg\'ebrique) 
de $X$-pr\'e\gro\  descend sur $\hat C$ et que $\hat C$ est un $X$-\gro (alg\'ebrique), avec 
$[c]^{-1} = [\bar c]$. 
La topologie sur $\hat C$ sera obtenue de la mani\`ere suivante~: consid\'erons
l'ensemble $\calk$ des paires $(K,T)$, o\`u 
\begin{itemize}
\item $K$ est un sous-groupe normal de $\Omega_{\hat C}$. Le quotient de $\hat C$
par le sous-\gro\ normal $N_K:=\{u K u^{-1}\mid u\in \hat C_*\}$ est alors un 
$X$-groupo\"\i de.
\item $T$ est une topologie sur $\hat C/N_K$ qui fait de $\hat C/N_K$ 
un \gro\ topologique
et telle que la projection $\tilde C \to \hat C/N_K$ soit continue. 
\end{itemize}
Les projections $\hat C \to \hat C/N_K$ ($(K,T)\in\calk$) forment un syst\`eme projectif de 
$X$-\gros\ au dessous de $\hat C$
(non-vide, car on peut prendre pour $T$ la topologie grossi\`ere). 
On munit $\hat C$ de la topologie de limite
projective pour ce syst\`eme de projections. On v\'erifie que $\hat C$ est alors
un $X$-groupo\"\i de, que $\tilde C \to \hat C$ est continue et que tout morphisme continu
de $\tilde C$ dans un $X$-\gro\ se factorise de fa\c con unique par l'un des quotient
$\hat C/N_K$. Comme $\tilde C$ est localement trivial num\'erisable, $\hat C$ l'est aussi.
Nous ignorons si la topologie ainsi obtenue sur $\hat C$ est la topologie quotient de celle de
$\tilde C$. 

Comme $X$ est s\'epar\'e, l'adh\'erence de $\{i_*\}$ est contenue dans $\Omega_C$ o\`u elle 
constitue un sous-groupe ferm\'e. En quotientant 
$\hat C$ par le sous-\gro\ normal engendr\'e par 
$\overline{\{i_*\}}$, on obtient, avec la topologie quotient, un $X$-\gro\ $C$ 
\cite[Th. 2.15, p. 38]{Ma}. Comme $\overline{\{i_*\}}$ est 
ferm\'e dans $\Omega_C$, le groupe $\Omega_C$ est s\'epar\'e. On en d\'eduit que
$C$ est s\'epar\'e par le lemme \ref{grosep}.
Il est aussi clair que que tout morphisme continu
de $\tilde C$ dans un $X$-\gro\ s\'epar\'e factorise de fa\c con unique par $\tilde C\to C$.

Soit $\tilde w : \tilde C\times E \to E$ une repr\'esentation de $\tilde C$ sur un $G$-fibr\'e 
principal $\xi$. Comme dans la proposition \ref{reptoto}, $\tilde w$ d\'etermine
un morphisme continu $\Psi_{\tilde w}$ de $\tilde C$ dans $EE:=EE(\xi)$ tel que 
$\Psi_{\tilde w}(\bar c) = \Psi_{\tilde w}(c)^{-1}$. Puisque $EE$ est localement
trivial et que $X$ et $G$ sont
s\'epar\'es, l'espace $EE$ est s\'epar\'e par le lemme \ref{grosep}. Le morphisme
$\Psi_{\tilde w}$ se factorise en un unique morphisme continu 
$\Psi_w:C\to EE$ correspondant, par la proposition \ref{reptoto}, \`a 
$w\in\calar$. \cqfd

\subsection{Groupo\"\i des de chemins} \label{SP:groch}

\begin{ccote}Chemins de Moore : \ \rm
Soit $X$ un espace topologique s\'epar\'e. 
Soit $\tilde C = \tilde {\bf Ch}(X)$  la cat\'egorie des {\it
chemins de Moore} dans $X$. Un chemin de Moore est un couple $(a,c)$ o\`u
$a\in\bbr_{\geq 0}$ et
$c : [0,a] \fl{}{} X$ est une application continue. 

La topologie sur $\tilde C$ est induite par l'application  
$(a,c) \mapsto c^{\sharp}$ de 
$\tilde C$ dans ${\rm map\,}([0,1],X)$, o\`u
$c^{\sharp}(t)=c(at)$, avec la CO-topologie sur ${\rm map\,}([0,1],X)$).
(Observons que $\tilde C$ n'est pas s\'epar\'e puisque tous les chemins constants
ont m\^eme image dans ${\rm map\,}([0,1],X)$).
Les applications
$\alpha,\beta : C\to X$ sont donn\'ees par $\alpha(a,c)=c(a)$ et $\beta
(a,c)=c(0)$. L'espace $C_x^y$ est donc l'ensemble des chemins allant de
$y$ \`a $x$ (cette malencontreuse inversion est due \`a la convention
usuelle de la r\`egle de composition des chemins).
La composition $(a,c)=(a_2,c_2)(a_1,c_1)$, lorsque $\beta
(a_1,c_1)=\alpha (a_2,c_2)$ est d\'efinie par $a:=a_2+a_1$ et 
 $$ c(t) = \cases{ c_2(t) & si $t \leq a_2$ \cr
                   c_1(t-a_2) & si $t \geq a_2$ \cr }$$
Cette composition est bien associative et l'unit\'e $i_x$
est donn\'ee par le chemin constant $[0,0] \fl{}{} \{x\}$. 
La d\'efinition de l'involution est donn\'ee par $\overline{(a,c)}:=(a,\bar c)$ o\`u
$\bar c(t):=c(a-t)$. On v\'erifie facilement que ${\bf Ch}(X)$ est un pr\'egroupo\"\i de. 
Il est localement trivial si et seulement si $X$ est connexe par arc et 
semi-localement contractile. 

Le \gro s\'epar\'e associ\'e \`a $\tilde {\bf Ch}(X)$ par la proposition \ref{predefgro}
sera not\'e ${\bf Ch}(X)$
et appel\'e le {\it \gro\ des chemins de Moore dans $X$}. 
\end{ccote}

\begin{ccote} Le \gro fondamental : \rm
Soit ${\bf Ch}_1(X)\subset {\bf Ch}(X)$ l'ensemble des classes de chemins
$(a,c)$ qui sont des lacets ($c(0)=c(a)$) et tels qu'il existe une homotopie 
$H:[0,a]\times [0,1]\to X$ telle que $H(0,s)=H(a,s)=c(0)$, $H(t,0)=c(t)$ et
$H(t,1)=c(0)$. Les \'el\'ements de ${\bf Ch}_1(X)$ forment un sous-$X$-\gro\ normal
totalement intransitif de ${\bf Ch}(X)$. L'espace quotient 
$\pi(X)$ h\'erite donc d'une structure de $X$-\gro\ et la projection
${\bf Ch}(X)\to\pi(X)$ est un morphisme continu \cite[Th. 2.15, p. 38]{Ma}.

Il est clair qu'alg\'ebriquement, $\pi(X)$ s'identifie au \gro\ fondamental 
de $X$ et $\Omega_{\pi(X)}$ au groupe fondamental $\pi_1(X,*)$ 
\cite[Ch. 1, \S\ 7]{Sp}. Cependant, $\pi(X)$ est ici muni d'une topologie.
Il est localement trivial si $X$ est connexe par arc et 
semi-localement simplement connexe \cite{Sp}. 
Si, de plus, $X$ est localement
connexe par arc, on peut montrer que $\Omega_{\pi(X)}=\Pi_1(X,*)$ est discret,
que le $\Omega_{\pi(X)}$-fibr\'e principal du th\'eor\`eme A est
le rev\^etement universel de $X$
et que la topologie sur $\pi(X)$ s'identifie \`a celle de \cite{BD}. 
Nous n'utiliserons pas ces r\'esultats. La proposition suivante est int\'eressante 
pour les $G$-fibr\'es principaux avec $G$ un groupe de Lie.
\end{ccote}

\begin{Proposition}\label{petitssgr}
Soit $G$ un groupe topologique admettant un voisinage de son \'el\'ement neutre 
qui ne contienne aucun sous-groupe non-trivial. Alors, toute repr\'esentation 
de ${\bf Ch}(X)$
sur un $G$-espace principal $\xi$ au dessus de $X$ se factorise par une repr\'esentation 
du \gro\ fondamental $\pi(X)$ de $X$.
\end{Proposition}

\preu Soit $w\in\calr ({\bf Ch}(X),\xi)$. On regarde $w$ comme un morphisme continu 
$w:{\bf Ch}(X)\to EE(\xi)$ par le lemme \ref{reptoto}. Il s'agit de montrer que
${\bf Ch}_1(X)\subset \ker w$.

Pour $x\in X$, d\'esignons par ${\rm Vois\,}_x$ l'ensemle des voisinages ouverts de $x$.
Observons que l'ensemble $\{\tilde {\bf Ch}(W)^x_x\mid W\in {\rm Vois\,}_x\}$ constitue
un syst\`eme fondamental de voisinages ouverts de $i_x$ dans le mono\"\i de
$\tilde {\bf Ch}(X)^x_x$.
En choisissant $\tilde x\in E_x$, on obtient une holonomie 
$h^w_x : \tilde {\bf Ch}(X)^x_x \to G$ en $x$
qui est un morphisme continu de mono\"\i des avec anti-involution. 
Comme il existe un voisinage de l'\'el\'ement neutre dans $G$ qui ne contient aucun sous-groupe
non-trivial et que chaque \'el\'ement de $\{\tilde {\bf Ch}(W)^x_x\mid W\in {\rm Vois\,}_x\}$
est un sous-mono\"\i de avec anti-involution de $\tilde {\bf Ch}(X)^x_x$, on en d\'eduit qu'il existe
$U_x\in {\rm Vois\,}_x$ tel que 
$\tilde {\bf Ch}(U_x)^x_x\subset \ker h^w_x$.
Soit ${\bf Ch}_{\rm loc}(X)$ le sous-\gro\ normal de ${\bf Ch}(X)$ engendr\'e par
l'image des $\tilde {\bf Ch}(U_x)^x_x$ pour tous les $x\in X$. Par ce qui pr\'ec\`ede, on a 
${\bf Ch}_{\rm loc}(X)\subset \ker w$.

Soit $\gamma\in {\bf Ch}_1(X)^x_x$ repr\'esent\'e par un lacet $c: [0,a]\to X$ en $x$.
Par d\'efinition de ${\bf Ch}_1(X)$, il existe une homotopie de lacets 
$H:[0,a]\times [0,1]\to X$ entre $c$ et le lacet constant. Par l'argument habituel du nombre
de Lebesgue, on peut d\'ecomposer $[0,a]\times [0,1]$ en petits rectangles $R_i$ ($1=1,\dots , N$)
tels que $H(R_i)\subset U_{x(i)}$. Il s'en suit facilement que 
$\gamma\in {\bf Ch}_{\rm loc}(X)\subset \ker w$, ce qui prouve que 
${\bf Ch}_1(X)\subset \ker w$.  \cqfd

\subsection{Chemins lisses par morceaux -- Connexions}\label{clipamconn}

Soit $X$ une vari\'et\'e diff\'erentiable $C^1$ paracompacte. On
d\'enotera par $\tilde {\bf D} = \tilde {\bf D}(X)$ l'espace des chemins
$C^1$ par morceau sur $X$. En tant qu'ensemble, $\tilde {\bf D}$ est le sous-pr\'e\gro\
de $\tilde {\bf Ch}X$ form\'e des chemins
qui sont des compositions de chemins $C^1$. 
La topologie, plus fine, s'obtient via l'application
$(a,c) \mapsto c^\sharp$ en topologisant 
l'ensemble ${\bf M}^\sharp$ des applications
$C^1$ par morceau de $[0,1]$ dans $X$. Pour cela, soit 
$$P = \{t_0=0 < t_1 < \cdots < t_k=1\}$$ 
un partage de $[0,1]$. Soit 
$${\bf M}_P^\sharp := \{c \in {\bf M}^\sharp \mid c\,_{|\,[t_i,t_{i+1}]} \in 
C^1([t_i,t_{i+1}],X) \} $$ 
o\`u $C^1([a,b],X)$ est l'espace des applications $C^1$ de $[a,b]$ dans $X$ muni de 
la topologie $C^1$. 
On a une application ${\bf M}_P^\sharp \to C^1([0,1],X)^k$ donn\'e par 
$$c \mapsto \big( c\,_{|\,[t_0,t_1]}^\sharp ,c\,_{|\,[t_1,t_2]}^\sharp , \dots
,c\,_{|\,[t_{k-1},t_k]}^\sharp \big)$$
en \'etendant la d\'efinition de $c^\sharp$ \`a une chemin $c :
[a,b]\fl{}{} X$ par $c^\sharp(t):=c(a+(b-a)t)$. L'application  $(a,c) \mapsto c^\sharp$
induit une topologie (non-s\'epar\'ee) sur ${\bf M}_P^\sharp$. 
Si $P'$ est un partage plus fin que $P$
(i.e. $P \subset P'$), on v\'erifie que l'inclusion naturelle 
${\bf M}_P^\sharp \subset {\bf M}_{P'}^\sharp$ est continue. La topologie sur
${\bf M}$ est, par d\'efinition,  celle de limite inductive des ${\bf M}_P^\sharp$
pour tous les partages de $[0,1]$.

Le \gro s\'epar\'e associ\'e \`a $\tilde {\bf D}(X)$ sera not\'e ${\bf D}(X)$
et appel\'e le {\it \gro\ des chemins lisses par morceaux dans $X$}. 
Si $X$ est connexe, alors ${\bf D}(X)$ est localement trivial num\'erisable. 
On a un morphisme \'evident de \gros topologiques de ${\bf D}(X)$ dans 
${\bf Ch}(X)$.

\begin{Proposition}\label{trans}
Soit $G$ un groupe de Lie et 
$\xi : E\hfl{p}{} X$ un $G$-fibr\'e principal diff\'erentiable $C^1$ 
au dessus d'une vari\'et\'e $X$.
Soit $A$ une connexion sur $\xi$ \cite[Ch. II]{KN}.
Alors, le transport parall\`ele associ\'e \`a $A$ d\'etermine une 
repr\'esentation de ${\bf D}(X)$ sur $\xi$.
\end{Proposition}

\preu Par la proposition \ref{predefgro}, il suffit de voir qu'une connexion $A$
d\'efinit une repr\'esentation $w_A : \tilde {\bf D}(X)\times E \to E$ de 
$\tilde {\bf D}(X)$, le point $w_A(c,z)$ \'etant le r\'esultat du transport
$A$-parall\`ele de $z$ au dessus de $c$. Les conditions 1. \`a 4. de la d\'efinition 
d'une repr\'esentation d\'ecoulent imm\'ediatement
des propri\'et\'es classiques du transport parall\`ele 
\cite[Ch. II, prop. 3.2. et 3.3]{KN}.
La seule chose \`a v\'erifier est que $w_A$ est continue en tout 
$(c,z)\in \tilde {\bf D}(X)\times E$. Vu la topologie sur $\tilde {\bf D}(X)$, il est 
suffisant de la faire pour $c: \in C^1([0,1],U)$ o\`u $U$ est un ouvert
de $X$ trivialisant pour $\xi$ et domaine d'une carte. 
On peut donc supposer que $E=U\times G$ ou $U$ est un ouvert d'un espace euclidien
et $z=(c(0),e)$.
Le relev\'e horizontal 
$\tilde c $ de $c$ partant de $z$ s'\'ecrit alors 
$\tilde c (t) = (c(t), g(t))$ o\`u $g\in C^1([0,1],G)$ avec $g(0)=e$.

Consid\'erons un voisinage ouvert $Q$ dans $E$ de $w_A(c,z)=\tilde c (1)= (c(1),g(1))$. 
Il s'agit de tourver un voisinage $T$ de $(c,z)$ dans $C^1([0,1],U)\times E$ 
tel que $w_A(T)\subset Q$. On peut supposer que $Q$ est 
de la forme $S\times (g(1)\cdot V)$ o\`u $V\in\calv_G$ o\`u $S$ est un ouvert de $U$. 
Soit $\varepsilon >0$ tel que la boule ouverte $B(0,\varepsilon)$ de centre $0$ et de rayon 
$\varepsilon$ dans $T_eG={\rm Lie\,}(G)$, muni de la m\'etrique de Killing, est envoy\'ee
diff\'eomorphiquement sur $W\in\calv_G$ avec $W\cdot W\subset V$.

Par \cite[Ch. II, prop. 1.1]{KN}, la connexion $A$ est donn\'ee par une $1$-forme
$\gamma_A\in\Omega^1(E,{\rm Lie\,}(G))$ et l'on a 
$\gamma_A(\dot{\tilde c}(t)) = 0$ pour tout $t$ puisque 
$\tilde c$ est horizontal. Par continuit\'e de $\gamma_A$, il existe
$\delta>0$ tel que $B(c(1),\delta)\subset S$ et
\begin{equation}\label{norme1}
\sup_{t\in [0,1]} \{\|c_1(t)-c(t)\|,\|\dot c_1(t)-\dot c(t)\|\}<\delta
\quad \Rightarrow \quad \|\gamma_A(\dot{\bar c_1}(t))\|<\varepsilon
\end{equation}
o\`u ${\bar c_1}(t):=(c_1(t),g(t))$. Soit $h\in C^1([0,1],G)$ la courbe telle
que $\hat c_1(t)\cdot h(t)$ soit horizontal. 
Par \cite[preuve du lemme p. 69]{KN}, la courbe $h$ satisfait $h(0)=e$ et
$\dot h(t) = T_e R_{h(t)}(\gamma_A(\dot{\bar c_1}(t))$, o\`u $R_a$ 
est la translation \`a droite $g\mapsto ag$. On
d\'eduit de \eqref{norme1} que $\ell (h)$ est $<\varepsilon$, o\`u
$\ell(h)$ est la longueur de $h$ pour 
la m\'etrique riemannienne sur $G$ obtenue par 
translations \`a droite de la m\'etrique de Killing.
L'exponentielle des rayons donnant des g\'eod\'esiques minimisantes pour
cette m\'etrique riemannienne, on en d\'eduit que $h(1)\in W$. 
L'ouvert $T:=T_\delta\times (B(c(0),\delta)\times W)$ de $C^1([0,1],U)\times G$, o\`u
$T_\delta$ est l'ouvert de $C^1([0,1],U)$ apparaissant dans \eqref{norme1}, contient
$(c,z)$ et satisfait $w_A(T)\subset Q$. \cqfd

Le langage des repr\'esentations de \gros\ permet de bien poser le probl\`eme suivant~:
quand est-ce que le transport parall\`ele associ\'e \`a une connexion sur un
fibr\'e diff\'erentiable s'\'etend aux chemins $C^0$~? La r\'eponse est la suivante~:

\begin{Proposition}\label{plate}
Soient $G$, 
$\xi$ et $A$ comme dans la proposition \ref{trans}.
Alors, la repr\'esentation $w_A\in\calr ({\bf D}(X),\xi)$ induite par 
le transport parall\`ele de $A$
s'\'etend en un repr\'esentation de ${\bf Ch}(X)$ sur $\xi$ si et seulement si 
$A$ est une connexion plate.
\end{Proposition}

\preu Par d\'efinition, $A$ est plate si et seulement si et seulement si 
$E$ est feuillet\'ee en vari\'et\'es horizontales \cite[II.9]{KN}. Il est clair que dans ce cas
le transport parall\`ele de $A$ s'\'etend en un repr\'esentation 
$\bar w_A\in\calr({\bf Ch}(X),\xi)$.
R\'eciproquement, si une telle extension existe, comme le groupe de Lie $G$ n'a pas 
de petits sous-groupes, la proposition \ref{petitssgr} assure que $\bar w_A$ se factorise
par le \gro\ fondamental $\Pi(X)$. Le th\'eor\`eme de r\'eduction classique \cite[II, th. 7.1]{KN}
montre qu'alors le fibr\'e $\xi_{\mid U} : p^{-1}(U)\to U$, au dessus de tout ouvert 
contractile $U$ de $X$, admet une r\'eduction de son groupe structural au groupe
trivial et que $A$ restreinte \`a $p^{-1}(U)$ est plate. \cqfd

\subsection{Th\'eorie de jauge sur graphes} \label{SP:lattice}

Rappelons qu'un {\it graphe} (non-orient\'e) $\Gamma$ 
consiste en une paire d'ensembles
$(S(\Gamma),A(\Gamma))$ (sommets et ar\^etes) avec deux applications 
$\alpha,\beta : A(\Gamma)\to S(\Gamma)$ 
et une involution $a\mapsto \bar a$ sur $A$ 
telle que $\alpha(\bar a) = \beta(a)$ et $\beta(\bar a)=\alpha(a)$.
Par exemple, pour $n\in\bbn$, le graphe $[n]$ se d\'efinit par
$S([n]):=\{0,1,\dots ,n\}$, $A([n]):=\{(i,j)\mid |i-j|=1\}$,
$\alpha(i,j):=j$, $\beta(i,j):=i$ et $\overline{(i,j)}:=(j,i)$.

Un {\it chemin} (de longueur $n$) dans $\Gamma$ est un morphisme de graphes
 de $[n]$ dans $\Gamma$. L'ensemble des chemins dans $\Gamma$ forment un 
$S(\Gamma)$-pr\'e\gro\ $\widetilde {C_\Gamma}$, muni
de la topologie discr\`ete. Les applications source et but, la composition et 
l'anti-involution sont d\'efinies comme pour les chemins de Moore 
dans un espace topologique (voir \S \ref{SP:groch}). Observons que 
$A(\Gamma)$ s'identifie naturellement au sous-ensemble de $\widetilde C_\Gamma$ 
form\'e des chemins de longueur 1.
Le $S(\Gamma)$-\gro associ\'e par la proposition \ref{predefgro}
sera not\'e ${\bf C}(\Gamma)$. Il est \'egalement discret et s'identifie au
\gro\ fondamental de la r\'ealisation g\'eom\'etrique $|\Gamma |$ de $\Gamma$. 
De m\^eme, $\Omega_{{\bf C}(\Gamma)}$ s'identifie au groupe fondamental 
$\pi_1(|\gamma |,*)$.
Nous supposerons que $|\Gamma |$ est connexe.

Ce qui s'appelle en anglais un "$G$-valued lattice gauge field" sur $\Gamma$
correspond \`a un $G$-fibr\'e $\xi$ sur $S(\Gamma)$ muni d'une repr\'esentation
$w\in\calr({\bf C}(\Gamma),\xi)$. Comme $S(\Gamma)$ est discret, 
le fibr\'e $\xi$ est trivial. Une trivialisation peut en \^etre obtenue \`a 
\`a l'aide de $w$ en choisissant un arbre maximal $T$ dans $\Gamma$.
En effet, $T$ donne une ${\bf C}(\Gamma)$-contraction sur tout $S(\Gamma)$.
En fixant une trivialisation
$E(\xi)=S(\Gamma)\times G$ de $\xi$, la repr\'esentation $w$ est d\'etermin\'ee par
la donn\'ee d'une application $w_1: A(\Gamma)\to G$ telle que $w_1(\bar a)=w_1(a)^{-1}$.
La formule reliant $w$ \`a $w_1$ est la suivante~:
\begin{equation}
w(a,(\alpha(a),g))=(\beta(a),w_1(a)\,g) \ , \ a\in A(\Gamma).
\end{equation}
Dans la litt\'erature sur le sujet, l'application 
$w_1$ est prise comme d\'efinition  d'un 
"$G$-valued lattice gauge field" (\cite[Ch. 7]{Cr}, \cite[\S\ 3]{PS}).

Si $\gamma$ est fini, le groupe $\pi_1(|\gamma |,*)$ est libre de rang 
$1-\chi(|\Gamma |)$, o\`u $\chi(|\Gamma |)$ est la caract\'eristique d'Euler de $|\Gamma |$.
Le th\'eor\`emes B et C donnent ainsi des hom\'eomorphismes
\begin{equation}
\calr ({\bf C}(\Gamma),\xi)/\calg_1 \approx \calr (\pi_1(|\gamma |,*),G) \approx
G^{1-\chi(|\Gamma |)}.
\end{equation}
et
\begin{equation}\label{gau2}
\calr ({\bf C}(\Gamma),\xi)/\calg \approx
G^{1-\chi(|\Gamma |)}\big/\hbox{\scriptsize conjugaison} .
\end{equation}

Sous certaines hypoth\`eses, la donn\'ee de $(\xi,w)$ d\'etermine un $G$-fibr\'e principal
$\xi_w$ sur $|\Gamma |$ qui n'est pas trivial (voir, par exemple, \cite{PS}). Par 
\eqref{gau2}, on obtient 
une partition de $G^{1-\chi(|\Gamma |)}\big/\hbox{\scriptsize conjugaison}$ en fonction 
des classes d'isomorphisme de $\xi_w$ qu'il serait int\'eressant d'\'etudier.

\subsection{Fibr\'es \'equivariants} \label{SP:fequi}

Soit $\Gamma$ un groupe topologique agissant \`a gauche sur $X$. 
Le graphe de l'action
$$C:=\{(y,\gamma,x)\in X\times\Gamma\times X\mid y=\gamma x\}$$
est un $X$-\gro\ par les donn\'ees suivantes~: 
$\alpha(y,\gamma,x):=x$, $\beta(y,\gamma,x):=y$, $i_x:=(x,e,x)$,
$(z,\gamma_2,y)(y,\gamma_1,x):=(z,\gamma_2\gamma_1,x)$ et 
$(y,\gamma,x)^{-1}:=(x,\gamma^{-1},y)$.
Le groupe $\Omega_{C}$ est banalement isomorphe au groupe d'isotropie 
$\Gamma^{\{*\}}$ de $*$.
Le $X$-\gro\ $C$ est localement trivial si l'action de $\Gamma$ sur $X$ est transitive et
si l'application $q:\Gamma\to X$ donn\'ee par $\gamma\mapsto \gamma\, *$ admet des sections
locales continues. Cette application sera alors un $\Gamma^{\{*\}}$ (le fibr\'e $\xi_C$ du th\'eor\`eme A).
On sait que cette situation se produit si, par exemple, 
$X$ est le quotient $\Gamma/\Gamma_0$ d'un groupe de Lie $\Gamma$ par un sous-groupe
ferm\'e $\Gamma_0$ \cite[\S\ 7.5]{St}.

Soit $\xi: E\hfl{p}{} X$ un $G$-fibr\'e principal sur $X$. 
Une repr\'esentation de $C$ sur $\xi$
est simplement une action \`a gauche de $\Gamma$ sur $E$ telle que la projection
$p:E\to X$ soit \'equivariante. On parle de $G$-fibr\'e principal $\Gamma$-\'equivariant
ou d'action de $\Gamma$ sur $\xi$.
L'espace $\calr (C,\xi)$ classe ces
$\Gamma$-actions sur $\xi$ \`a conjugaison par une transformation de jauge pr\`es.
Par le th\'eor\`eme B, $\calr (C,\xi) \approx \calr (\Gamma^{\{*\}},G)_\xi$.
Remarquons que $B\times_{\Gamma^{\{*\}}} E\Gamma \approx B\Gamma^{\{*\}}$; dans le cas o\`u 
$G$ est ab\'elien, on retrouve ainsi le th\'eor\`eme A de [LMS]. Les relations avec d'autres approches
de fibr\'es \'equivariants, comme par exemple [BH], restent \`a \'etudier.

Voici quelques exemples : 

\begin{ccote}
$\Gamma = SO(3)$ agissant sur $X=S^2$ et $G=S^1$. \rm
On a $\Omega_{C}=S^1$ et le fibr\'e principal 
$\xi_{C}$ du th\'eor\`eme A est le fibr\'e tangent unitaire \`a $S^2$,
dont la classe d'Euler est 2. Par le th\'eor\`eme d'existence,
un $S^1$-fibr\'e principal sur $S^2$ admettra une $SO(3)$-action 
si et seulement si sa classe d'Euler est paire. Dans ce cas, 
l'espace $\calr (S^1,S^1)$ \'etant discret
(hom\'eomorphe \`a $\bbz$ par le degr\'e), le th\'eor\`eme B implique qu'il y a exactement
une classe d'action de $SO(3)$ sur $\xi$ 
\`a conjugaison par une transformation de jauge pr\`es. Observons qu'il n'y
a aucun choix possible pour l'action sur la fibre au dessus de $*$~; cette action est d\'etermin\'ee par
$\xi$. 
\end{ccote}

\begin{ccote}
$\Gamma = SU(2)$ agissant sur $X=S^2$ et $G=S^1$. \rm
Le fibr\'e $\xi_{C}$ du th\'eor\`eme A est alors le fibr\'e de Hopf $S^3\to S^2$. 
On en d\'eduit que tout $S^1$-fibr\'e principal sur $S^2$ admet une $SU(2)$-action 
unique \`a conjugaison par une transformation de jauge pr\`es. 
\end{ccote}

\begin{ccote}
$\Gamma = SU(2)$ agissant sur $X=S^2$ et $G=SO(3)$. \rm
Il y a deux $SO(3)$-fibr\'es principaux $\xi_0$ et $\xi_1$ sur $S^2$, $\xi_0$ \'etant le fibr\'e
trivial. Tout deux associ\'es au fibr\'e de Hopf, ils admettent des $SU(2)$-actions mais seul
$\xi_0$ admet des $SO(3)$-actions. 
Un homomorphisme continu de $S^1$ dans $SO(3)$ est diff\'erentiable, 
donc un \'el\'ement de $\calr(S^1,SO(3))$ est un sous-groupe 
\`a un param\`etre dont l'image est un tore maximal
de $SO(3)$. On en d\'eduit que si
l'on identifie l'alg\`ebre de Lie $so(3)$ \`a $\bbr^3$ (quaternions purs), l'espace
$\calr(S^1,SO(3))$ est la r\'eunion des sph\`eres de rayons entiers $n=0,1,2,\dots$.
Donc, $\calr(C,\xi_0)/\calg_1$ est hom\'eomorphe \`a la r\'eunion des 2-sph\`eres de rayons $2n$
($n=0,1,2,\dots$) et $\calr(C,\xi_1)/\calg_1$ \`a celles de rayon $2n+1$.
Quant \`a $\calr(C,\xi_i)/\calg$, ils sont tout deux discrets d\'enombrables.
\end{ccote}

\vskip .3 truecm\goodbreak


\vskip .5 truecm\small
\noindent \parbox[t]{6 truecm}{Jean-Claude HAUSMANN\\
Math\'ematiques-Universit\'e\\ B.P. 240, \\
 CH-1211 Gen\`eve
 24, Suisse\\ hausmann@math.unige.ch} \ \hfill \hfill \

\end{document}